%% file: thesis.tex
\documentclass[12pt,a4paper,BCOR=10mm]{scrbook} 

\usepackage{scrhack} 
\usepackage[utf8]{inputenc}
\usepackage[inner=3cm,outer=2cm,top=3cm,bottom=4cm]{geometry}
\usepackage[T1]{fontenc}
\usepackage{lmodern}
\usepackage{amsthm,amsmath,amssymb, mathtools}
\usepackage{bbm} 
\usepackage[autostyle=true]{csquotes} 
\usepackage{color}
\usepackage{thmtools}
\usepackage[english,ngerman]{babel} 
\usepackage{url}
\usepackage[backend=biber,citestyle=ieee,bibstyle=ieee]{biblatex}
\usepackage{booktabs}
\usepackage{pgfplots}
\usepackage{pgfplotstable}
\usepackage{caption}
\usepackage{subcaption}     
\usepackage{algorithm} 

\usepackage{chngcntr}
\counterwithout{footnote}{chapter}

\usepackage{microtype} 

\usepackage[titletoc]{appendix}
\usepackage[noend]{algpseudocode} 
\usepackage[hidelinks,hyperfootnotes=false,bookmarksopen=true]{hyperref}

\usepackage{imakeidx}
\makeindex[intoc]

\usepackage{setspace}
\pgfplotsset{
	compat=1.13,
}

\usetikzlibrary{external,calc,arrows,automata,positioning,intersections}

\newtheorem{theorem}{Theorem}[chapter]

\newtheorem{lemma}[theorem]{Lemma}
\theoremstyle{definition}
\newtheorem{definition}[theorem]{Definition}
\newtheorem{proposition}[theorem]{Proposition}
\newtheorem{observation}[theorem]{Observation}
\newtheoremstyle{theostyle}
	{5pt} 
	{5pt} 
	{} 
	{} 
	{\bfseries} 
	{.} 
	{.5em} 
	{} 
\theoremstyle{theostyle}
\declaretheorem[style=theostyle,sibling=definition]{example}

\newcommand*{\C}{\mathbb{C}}
\newcommand*{\R}{\mathbb{R}}
\newcommand*{\N}{\mathbb{N}}
\newcommand*{\eps}{\varepsilon}
\newcommand*{\iu}{{i\mkern1mu}}

\newcommand*{\ind}{\mathbbm{1}_{(-1,1)}}

\DeclareMathOperator{\Ima}{Ima}

\DeclareMathOperator{\am}{argmin}

\newcommand*\argmin[2]{\underset{#1}{\am} \ {#2}}
\newcommand{\set}[1]{\{\,#1\,\}}
\DeclarePairedDelimiter\abs{\lvert}{\rvert}
\newcommand*{\abbs}[1]{\left\lvert#1\right\rvert}

\newcommand*\colvecalt[2]{\begin{pmatrix}#1 \\ #2\end{pmatrix}}
\newcommand*\conj[1]{\overline{#1}}
\newcommand*\ctrans[1]{#1^{H}}
\newcommand*\trans[1]{#1^{T}}
\newcommand*\ratcp[2]{#1^{#2}}

\newcommand\litem[1]{\item{\bfseries #1.}}

\mathchardef\ordinarycolon\mathcode`\:
\mathcode`\:=\string"8000
\begingroup \catcode`\:=\active
\gdef:{\mathrel{\mathop\ordinarycolon}}
\endgroup

\newcommand{\expect}{\operatorname{E}\expectarg}
\DeclarePairedDelimiterX{\expectarg}[1]{[}{]}{%
	\ifnum\currentgrouptype=16 \else\begingroup\fi
	\activatebar#1
	\ifnum\currentgrouptype=16 \else\endgroup\fi
}
\newcommand{\innermid}{\nonscript\;\delimsize\vert\nonscript\;}
\newcommand{\activatebar}{%
	\begingroup\lccode`\~=`\|
	\lowercase{\endgroup\let~}\innermid 
	\mathcode`|=\string"8000
}

\newcommand{\tool}[1]{\textsf{#1}}
\newcommand*{\RFF}{\tool{RFF}}
\newcommand*{\wina}{\textsc{Winkelmann} and \textsc{Di Napoli}}

\newcommand*{\HIEP}{\tool{HIEP}}
\newcommand*{\CHIEP}{\tool{CHIEP}}
\newcommand*{\BFGS}{\tool{BFGS}}
\newcommand*{\LBFGSB}{\tool{L-BFGS-B}}
\newcommand*{\FEAST}{\tool{FEAST}}

\definecolor{LSOPC}{RGB}{77,175,74} 
\definecolor{GAUSSC}{RGB}{55,126,184} 
\definecolor{ZOLOC}{RGB}{228,26,28} 
\definecolor{LSOPC2}{RGB}{152,78,163} 
\definecolor{LSOPC3}{RGB}{255,127,0} 
\definecolor{XAXIS}{RGB}{100,100,100}

\def\LSOP/{\tool{SLiSe}} 

\usepackage{quotchap}
\newcommand{\preamble}[3]{\begin{savequote}[0.86\linewidth]
		\enquote{#3}
		\qauthor{--- \textsc{#1}, #2}
\end{savequote}}
\newcommand{\preambletrans}[4]{\begin{savequote}[0.86\linewidth]
		\enquote{#3} \\ \\		
		\enquote{#4}
		\qauthor{--- \textsc{#1}, #2}
\end{savequote}}

\bibliography{references}

\begin{document}
	\pagenumbering{gobble} 
	\include{parts/title}

	\newpage~\cleardoublepage
	\include{parts/ackstatement}		
	\include{parts/abstract}
	
	\microtypesetup{protrusion=false}
	\setstretch{1.41}
	\tableofcontents
	\onehalfspacing
	\cleardoublepage
	\microtypesetup{protrusion=true}
	
	\pagenumbering{arabic}	
	\include{parts/introduction}
	\include{parts/eigenvalues}	
	\include{parts/RFF}
	\include{parts/lossfunction}
	\include{parts/constraints}
	\include{parts/experiments}	
	\include{parts/conclusion}
	
	\pagenumbering{gobble} 
	\microtypesetup{protrusion=false}	
	\begin{appendices}
		\appendixpage
		\pagenumbering{roman}
		\include{parts/appendix}
	\end{appendices}
	\pagenumbering{gobble}

	\clearpage 
	\singlespacing 
	\addcontentsline{toc}{chapter}{Bibliography}
	\printbibliography
	\printindex
\end{document}

%% file: parts/title.tex
\title{
	\makebox[\textwidth][r]{\includegraphics[width=0.5\textwidth]{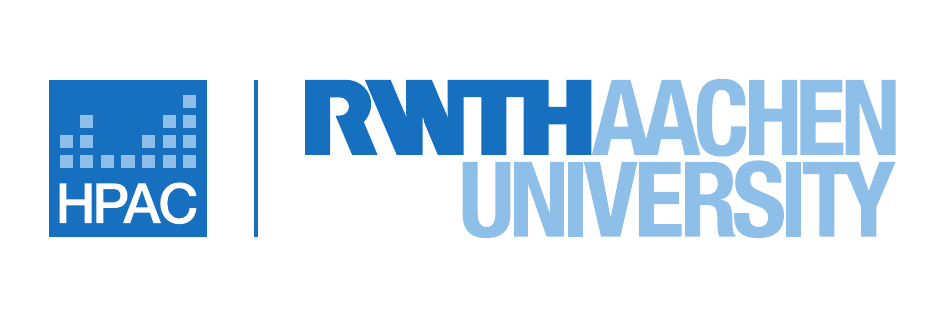}} \\  \vspace{-1cm} {\small \normalfont Diese Arbeit wurde vorgelegt am Lehrstuhl für High-Performance and Automatic Computing.}
	\vspace{4cm} \\
	{\large Bachelorarbeit}
	\vspace{1cm}
	\hrule \vspace{0.5cm}
	\normalfont \bfseries Constrained Optimisation of \\Rational Functions for \\ Accelerating Subspace Iteration \vspace{0.5cm}
	\hrule
}
\author{
	Konrad A. S. Kollnig \\
	\vspace{2cm}
}
\date{
	\selectlanguage{ngerman}
	\large Aachen, den 21. August 2017
}

\titlehead{
	\centering
}

\publishers{\begin{tabular}{rl}
	Erstprüfer:& Prof.\ Paolo Bientinesi, Ph.D. \\
	Zweitprüfer:& Prof.\ Georg May, Ph.D. \\
	Betreuer:& Jan Winkelmann, M.Sc.
	\end{tabular}
}

\begin{titlepage}
	\newgeometry{left=2.3cm, right=2.3cm, bottom=2cm,top=-2cm}
	\maketitle
	\restoregeometry
	\vspace*{10cm}
	\thispagestyle{empty}
	\clearpage
\end{titlepage}

%% file: parts/ackstatement.tex
\thispagestyle{empty}
\selectlanguage{ngerman}
\textbf{\sffamily \LARGE Danksagungen} 

Meine Zeit am \tool{HPAC}-Lehrstuhl habe ich sehr genossen. Die Vernetzung in alle Welt stellte eine Bereicherung mit neuen Eindrücken für mich dar. Dies war neben meiner Arbeit an dieser Ausarbeitung eine willkommene Abwechlung, um den Blick über den Tellerrand dieses Dokuments hinaus zu wagen.

Mein Dank gilt in erster Linie meinem Betreuer, \textsc{Jan Winkelmann}. Von ihm stammt das Gros der hier ausgearbeiten Ideen. Ohne seine fortwährende Unterstützung wäre diese Arbeit nicht möglich gewesen. Mit größter Sorgfalt stand er mir Antwort auf alle meine Fragen gleich welcher Stunde.


Großer Dank gilt auch \textsc{Prof. Paolo Bientinesi}, dem Leiter des \tool{HPAC}-Lehrstuhls. Seine Expertise war ein wertvolles Quell an Wissen. Geschätzt habe ich die Genauigkeit, mit der er die notwendigen Arbeitsschritte zu jeder Zeit herausstellen konnte.
\selectlanguage{english}

\vfill

\selectlanguage{ngerman}
{\LARGE \sffamily \textbf{Erklärung}}

Hiermit versichere ich, dass ich die vorliegende Arbeit selbständig verfasst, und keine anderen als die angegebenen Quellen und Hilfsmittel verwendet, sowie Zitate kenntlich gemacht habe.

\vspace{1cm}
	
Aachen, den 21. August 2017
	
\vspace{2cm}
	
\begin{tabular}{c}
	\hline
	$\qquad$ Konrad A. S. Kollnig $\qquad$ 
\end{tabular}

\vspace{2cm}
\selectlanguage{english}
\clearpage
\textbf{\sffamily \LARGE Acknowledgements}

My time at the \tool{HPAC} group has been a pleasure. I have enjoyed the international atmosphere and the exchange with researchers from all over the world. Besides my work on my thesis, this has been an enriching glimpse over the rim of the tea cup.

Foremost, I thank my supervisor, \textsc{Jan Winkelmann}. The majority of ideas discussed in this thesis are his achievement. Without his ongoing support, the creating of this document would have not been possible. Paying attention to the detail, he answered my questions even late at night.


Furthermore, I thank \textsc{Prof.\ Paolo Bientinesi}, the head of \tool{HPAC}. His expertise was a valuable source of knowledge. He could tell accurately, what to focus on during the work on this document.
\selectlanguage{english}

%% file: parts/abstract.tex
\cleardoublepage

\thispagestyle{empty}

\vspace*{\fill}

\begin{center}
	\textbf{\sffamily \LARGE Abstract} 
\end{center}
\hrule
Earlier this decade, the so-called \FEAST{} algorithm was released for computing the eigenvalues of a matrix in a given interval. Previously, rational filter functions have been examined as a parameter of \FEAST{}. In this thesis, we expand on existing work with the following contributions: (i) Obtaining well-performing rational filter functions via standard minimisation algorithms, (ii) Obtaining constrained rational filter functions efficiently, and (iii) Improving existing rational filter functions algorithmically. Using our new rational filter functions, \FEAST{} requires up to one quarter fewer iterations on average compared to state-of-art rational filter functions.
\vspace*{\fill}

\cleardoublepage

%% file: parts/introduction.tex
\preamble{E. W. Dijkstra}{EWD 1036}{The problems of the real world are primarily those you are left with when you refuse to apply their effective solutions.}
\chapter{Introduction}
\section{Rational filter functions for FEAST}
\label{sec:intro}
Eigenvalue problems are a fundamental subject of computational linear algebra. The \FEAST{}\index{FEAST} algorithm solves a particular case of eigenvalue problems, i.e.\ \textsc{Hermitian} interior eigenvalue problems \cite{feast}. By construction, \FEAST{} offers a natural support for parallelism to make use of state-of-art computer architectures. Recent discussion of the underlying \textit{parameters} of \FEAST{} has led to significant improvements in performance \cite{guettel,eigencount,feast_conv,jan,linear}. To achieve better performance, we examine \textit{rational filter functions}\index{Rational!Filter function|see{RFF}} (\RFF{}s\index{RFF}), a parameter\footnote{Other rational functions are accepted as well, but we want to build upon the currently used \RFF{}s.} of \FEAST{}, in this thesis.

An \RFF{} is a rational function, defined as
\begin{equation*}
	\ratcp{r}{\beta,w}(x) := \sum_{i=1}^{m} \frac{\beta_i}{x - w_i} + \frac{\conj{\beta_i}}{x - \conj{w_i}} - \frac{\beta_i}{x + w_i} - \frac{\conj{\beta_i}}{x + \conj{w_i}}, \quad \text{for} \ x \in \R,
\end{equation*}
where $\beta \in \C^m$, $w \in (\C \setminus \R)^m$, for some $m \in \N$.

The problem is finding appropriate \RFF{}s $\ratcp{r}{\beta,w}$, that perform well in \FEAST{} for a variety of eigenvalue problems. This is essential for good average performance of \FEAST{}.

Recently, \wina{} have introduced a new approach, called \textit{SLiSe minimisation}, to obtain well-performing \RFF{}s \cite{jan}. This approach aims at improving existing state-of-art \RFF{}s, such as \textsc{Gauss-Legendre} \cite{feast} or \textsc{Zolotarev} \cite{guettel}, by solving \textit{non-linear minimisation problems}\index{Minimisation problem}. Non-linear means, that a minimisation problem cannot be solved directly. In other words, non-linear problems are computationally difficult to solve.

An integral parameter of \LSOP/ are \textit{weight functions}\index{Weight function}, i.e.\ non-negative step functions. By weight functions, arbitrary degrees of freedom are offered in \LSOP/ minimisation to reflect and overcome the diversity of eigenvalue problems. It has been discussed, that the choice of weight function is crucial for obtaining well-performing \RFF{}s via \LSOP/. However, by a judicious selection of weight functions, a variety of \RFF{}s $\ratcp{r}{\beta,w}$ have been computed via \LSOP/ minimisation, outperforming the state-of-art \RFF{}s used in \FEAST{} \cite{jan}.

In summary, two aspects of \LSOP/ are of utmost importance: Solving the computationally difficult \LSOP/ minimisation problems efficiently and choosing suitable weight functions in \LSOP/. This currently results in the following drawbacks.
\begin{enumerate}
	\litem{Incompatibility with standard minimisation algorithms} The non-linear \LSOP/ minimisation problems require complex variables. This approach is problematic, as standard minimisation algorithms are incompatible with complex variables. Standard minimisation algorithms potentially yield better \RFF{}s, while decreasing the time required for solving the non-linear minimisation problems.
	
	\litem{Impracticable performance of box-constraints} To improve the performance of \RFF{}s in \FEAST{}, \wina{} suggested the usage of \textit{box-constraints}\index{Box-constraint} in \LSOP/ minimisation \cite{jan}. While they demonstrated the potential of this method, the performance of their implementation is not yet suitable for practice.
	
	\litem{Choice of suitable weight functions} The choice of suitable weight functions as part of \LSOP/ is fundamental for computing well-performing \RFF{}s. In the past, weight functions have been chosen empirically. If we had a criterion for estimating the performance of weight functions, we could algorithmically compare a multitude of weight functions to obtain better \RFF{}s via \LSOP/ minimisation.
\end{enumerate}

Our contributions solve these three issues of \LSOP/. We present them in the following.

\section{Contributions}
In this thesis, we solve the three issues of \LSOP/ introduced previously.
\begin{enumerate}
	\litem{Embedding into standard minimisation algorithms} We perform the embedding of \LSOP/ into standard minimisation algorithms in Chapter~\ref{sec:real}. As a prominent example of a minimisation algorithm, we study the \BFGS{} algorithm, which provides a good performance in our case.
	
	\litem{Improving the performance of box-constraints} We have improved the performance of solving the box-constrained \LSOP/ minimisation by the \LBFGSB{} algorithm. Compared to previous approaches, our approach works well for different \LSOP/ minimisation problems, requires fewer iterations in general and yields better \RFF{}s. We outline our approach and further improvements of \RFF{}s in Chapter~\ref{chp:constraints}.
	
	\litem{Algorithmic choice of weight functions}
	We propose a method to obtain well-performing weight functions for \LSOP/ algorithmically instead of empirically. This has enabled us to compute better \LSOP/ \RFF{}s. Our procedure and an analysis are presented in Chapter~\ref{chp:tests}. All our resulting \RFF{}s can be found in Appendix~\ref{chp:resulting_rffs}. They can be used in \FEAST{} directly. Our new weight functions to improve \LSOP/ are denoted in Appendix~\ref{chp:weight_functions}.
\end{enumerate}
\section{Structure}
In Chapters~$2-3$, we revise the required background briefly: In Chapter~\ref{chp:subspace}, we introduce the \FEAST{} algorithm. Notably, we focus on the convergence rate of \FEAST{} in Section~\ref{sec:feast_conv}. Reducing this convergence rate is our aim in this thesis. In Chapter~\ref{chp:rff}, we introduce existing \RFF{}s for further improvement in the following chapters.

In Chapters~$4-6$, we present our contributions: In Chapter~\ref{sec:real}, we embed \LSOP/ into standard minimisation algorithms. In Chapter~\ref{chp:constraints}, we examine constraints on \RFF{}s. In Chapter~\ref{chp:tests}, we obtain new weight functions for \LSOP/ algorithmically.

In Appendix A, all the discussed \RFF{}s can be found. In Appendix B, all the considered weight functions are denoted. In Appendix C, the \BFGS{} minimisation algorithm is discussed extensively. In Appendix D, we provide proofs of mentioned theoretical details.

\section{Preliminaries}

The following notation will be used throughout this thesis:
\begin{itemize}
	\item We denote the natural numbers excluding $0$ as $\N$.
	\item The complex numbers\index{Complex numbers} are defined as
	\begin{equation*}
	\C := \set{ a + b \cdot \iu \mid a,b \in \R},
	\end{equation*}
	where the symbol $\iu$ is called the \textit{imaginary unit}. For any complex number $z = a + b \cdot \iu$, $a,b \in \R$, we call $\Re(z) := a$ the \textit{real part} of $z$ and $\Im(z) := b$ the \textit{imaginary part} of $z$. Furthermore, we define the \textit{complex conjugate} as $\conj{z} := a - b \cdot \iu$ and the \textit{absolute value} as $\abs{z} := \sqrt{a^2 + b^2}$.
	\item For a given matrix $A$, we denote the \textit{conjugate transpose}\index{Conjugate transpose} of $A$ as $A^H = \conj{A^T}$. The conjugate of a matrix is defined as the conjugate of its elements.
	\item For a set $M$ and a function $f: M \rightarrow \R$, the notation
	\begin{equation*}
	\argmin{x \in M}{f(x)}
	\end{equation*}
	indicates, that we want to find an $x \in M$ that $f$ attains its minimum for. This is called the \textit{minimisation problem}\index{Minimisation problem} or \textit{optimisation problem}\index{Optimisation problem|see{Minimisation problem}} of $f$.
	\item When we first introduce new terminology, we use \textit{italics}. We explain terminology when first used.
	\item When we refer to equations, we use braces: \enquote{The \BFGS{} update is given by \eqref{eqn:bfgs}.}
	\item For a given function $f:A \rightarrow B$, we denote the \textit{image of $f$} as
	\begin{equation*}
		\Ima(f) := \set{f(x) \mid x \in A} \subseteq B.
	\end{equation*}
\end{itemize}

%% file: parts/eigenvalues.tex
\preambletrans{D. Hilbert}{Grundzüge einer allgeminen Theorie der linaren Integralrechnungen}{Dieser Erfolg ist wesentlich durch den Umstand bedingt, daß ich nicht, wie es bisher geschah, in erster Linie auf den Beweis für die Existenz der \textit{Eigenwerte} ausgehe.}{This success is mainly determined by the circumstance, that I do not, as happened previously, focus on the proof of the existence of \textit{eigenvalues} in the first place.}
\chapter[FEAST algorithm, subspace iteration with rational filter functions]{FEAST, subspace iteration with RFFs}
\label{chp:subspace}

To study the \FEAST{} algorithm, we first examine in Section~\ref{sec:hiep} the problem it solves, i.e. the \textsc{Hermitian} interior eigenvalue problem. For further improvement, we introduce the \FEAST{} algorithm in Section~\ref{sec:feast}. In Section~\ref{sec:feast_conv}, we discuss the convergence rate of \FEAST{}. This is what existing rational filter functions try to reduce for a variety of eigenvalue problems and what we will further dissect in this thesis.

\section{Problem}
\label{sec:hiep}
For the rest of this thesis, we assume a \textsc{Hermitian} matrix $A = A^H\in \C^{n \times n}$, for some $n \in \N$. The \textsc{Hermitian} \textit{eigenvalue problem}\index{Eigenvalue problem} is then identified by the equation
\begin{equation}
	\label{eqn:eigen}
	A v = \lambda v,
\end{equation}
for coefficients $\lambda \in \C$ and associated vectors $v \in \C^n \setminus \set{0}$ and some $n \in \N$.

We will call a solution $(\lambda,v)$ to this problem an \textit{eigenpair}\index{Eigenpair} of $A$. The coefficient $\lambda$ is known as an \textit{eigenvalue}\index{Eigenvalue} of $A$, the vector $v$ as an \textit{eigenvector}\index{Eigenvector} of $A$. It can be shown that a matrix $A$ has only finitely many different eigenvalues. Thus, it is possible to compute the eigenvalues for a given matrix $A$. It is sufficient to compute one corresponding eigenvector each.

In this thesis, we study the \textsc{Hermitian} \textit{interior} eigenvalue problem (\HIEP{})\index{Eigenvalue problem!Hermitian interior|see{HIEP}}\index{HIEP}. Interior means, that the eigenvalues lie in a bounded interval: For a given proper interval $(a,b)$, the \HIEP{} is computing the eigenpairs $(\lambda,v)$ of $A$ such that $\lambda \in (a,b)$, where $a,b\in\R$ fixed. We will refer to such an interval as the \textit{search interval} of the \HIEP{}. We focus on search intervals where $\lambda_\text{min} \le a,b \le \lambda_\text{max}$, where $\lambda_\text{min}$ denotes the minimum and $\lambda_\text{max}$ the maximum eigenvalue of $A$. It can be shown that all the eigenvalues of a \textsc{Hermitian} matrix are real-valued. This justifies searching for eigenvalues in a particular \textit{real} search interval only.

Note that any \HIEP{} can be reduced to a canonical form as follows. The search interval of a \HIEP{} can be phrased as ${(m-r,m+r)}$, where $m \in \R$ denotes the midpoint and $r \in (0,\infty)$ the radius of the given search interval. We then solve the \HIEP{} for the \textsc{Hermitian} matrix $A'$ given by
\begin{equation}
\label{eqn:reduce_filter_interval}
A' := (A - m \cdot I)/r
\end{equation}
instead of $A$ with search interval $(-1,1)$, where $I$ denotes the matching identity matrix.

This is shifting and scaling of the eigenvalues of $A$ in $(m-r,m+r)$ to $(-1,1)$. The eigenvectors of both the matrices $A$ and $A'$ remain unchanged. We will refer to a \HIEP{} with search interval $(-1,1)$ as \textit{Canonical \HIEP{}} (\CHIEP{})\index{CHIEP|see{HIEP}}. In the following, it is sufficient to focus on \CHIEP{}s only.

\section{Algorithm}
\label{sec:feast}

To introduce the \FEAST{} algorithm for solving \CHIEP{}s\footnote{\FEAST{} can also solve \textit{generalised} \HIEP{}s. We will only study the standard case. However, this entire document holds for the generalised case as well.}, we proceed as follows: First, we discuss the inputs of the algorithm. Then, we express the algorithm in terms of pseudo-code and point out the most important details for our analysis. Afterwards, we discuss implementation details for further understanding.

In addition to a \textsc{Hermitian} matrix, the \FEAST{} algorithm requires two main inputs: a rational filter function and an upper bound on the number of eigenvalues in the search interval. These inputs are introduced in the following.

\begin{enumerate}
	\litem{Rational filter function} For our analysis, the most important input of \FEAST{} is the rational filter function (\RFF{})\index{RFF}. An \RFF{} can be used by \FEAST{} to transform a \CHIEP{} to a \textit{reduced} (unbounded) eigenvalue problem. Reduced means that this eigenvalue problem contains fewer eigenvalues than the original matrix $A$. An \RFF{} is defined as
	\begin{equation}
	\label{eqn:rff}
	\ratcp{r}{\beta,w} := r(x) := \sum_{i=1}^{m} \frac{\beta_i}{x - w_i} + \frac{\conj{\beta_i}}{x - \conj{w_i}} - \frac{\beta_i}{x + w_i} - \frac{\conj{\beta_i}}{x + \conj{w_i}},
	\end{equation}
	for some $\beta \in \C^m$, $w \in (\C \setminus \R)^m$ and $m \in \N$.
	
	By assuming $w$ to be non-real, we assure that any \RFF{} is well-defined. The \RFF{} is said to have \textit{degree} $4m$. The values $\pm w_i, \ \pm \conj{w_i}$ are called \textit{poles}\index{Pole}\index{RFF!Pole} of $r$.
	
	An \RFF{} appears to map to complex values. However, from its definition, we note that $r(x) = \conj{r(x)}$, for all $x \in \R$. There exists no $z \in \C \setminus \R$ with $z = \conj{z}$. So, an \RFF{} is indeed \textit{real-valued}\index{RFF!Properties}. Furthermore, from its definition, we have $r(x) = r(-x)$, for all $x \in \R$. This means that \RFF{}s are \textit{symmetric with respect to the $y$-axis}\index{RFF!Symmetry}. These two properties are important, as shown in the analysis of the convergence rate in the following Section~\ref{sec:feast_conv}. More properties were proved in \cite{jan}.
	
	In the \FEAST{} algorithm, an \RFF{} is interpreted as a \textit{matrix function}. This is induced by a given \RFF{} $r$ as per
	\begin{equation}
	\label{eqn:rff_matrix}
	r(A) = \sum_{i=1}^{m} \beta_i (A - w_i I)^{-1} + \conj{\beta_i} (A - \conj{w_i} I)^{-1} - \beta_i (A + w_i I)^{-1} - \conj{\beta_i} (A + \conj{w_i} I)^{-1},
	\end{equation}
	where $I$ denotes the matching identity matrix.
	
	This matrix function is only evaluated when multiplied with another vector $v \in \C^n \setminus \set{0}$. Then, by computing the product $r(A)v$, $4m$ linear systems of the form
	\begin{equation}
	\label{eqn:feast_rff}
	\alpha (A - z I)^{-1} v,
	\end{equation}
	have to be solved, where $\alpha \in \C$ and $z \in \C \setminus \R$.
	
	The synergy of \FEAST{} and \RFF{} is important for the algorithm to function. However, for this thesis, this is out of scope, because we focus on the performance \RFF{}s in \FEAST{}. This performance is solely determined by the convergence rate, discussed in the following section. For more details on the usage of \RFF{}s in \FEAST{}, see  \cite{feast_conv}.
	\litem{Upper bound on eigenvalue number} As an additional input, \FEAST{} requires an upper bound $N \in \N$ on the number of eigenvalues in $(-1,1)$, the search interval of the \CHIEP{}. The algorithm requires knowledge on how many (interior) eigenpairs to search for. There exists an efficient technique to compute such an upper bound $N$ as part of the \FEAST{} algorithm, see \cite{eigencount}. Thus, in practice, no such prior knowledge of the \CHIEP{} is required.
\end{enumerate}

The core of the \FEAST{} algorithm is given as follows. It can also be called \textit{subspace iteration with rational filter functions}\index{Rational!Subspace Iteration}, as it generalises the well-known \textit{subspace iteration}\index{Subspace iteration} algorithm using \RFF{}s. For more details, confer \cite[ch. 5]{saad}.

\begin{algorithm}[H]
	\caption{(Subspace iteration with rational filter functions \cite{saad,feast,guettel,feast_subspace}).}
	\label{algo:subspace_iteration}
	\textbf{Input:} \textsc{Hermitian} matrix $A \in \C^{n \times n}$, an upper bound $N \in \N$ on the number of eigenvalues of $A$ in $(-1,1)$ and an \RFF{} $r$.\\
	\textbf{Output:} All eigenpairs $(\lambda,v)$ of $A$, where $\lambda \in (-1,1)$.
	\begin{algorithmic}[1]
		\Function{RffSubspaceIteration}{$A, N, r$}		
		\State $V \gets \Call{Random}{\C^{n \times N}}$ \Comment{Initial eigenvector approximation.}
		\Repeat
		\State $X \gets r(A) V$ \label{line:rff}
		\State $Q \gets \Call{Orthogonalise}{X}$ \label{line:ortho}
		\State $B \gets Q^H A Q$
		\State $(\Lambda,Y) \gets \Call{Solve}{B Y = B \Lambda}$ \Comment{Reduced eigenvalue problem.} \label{line:reduced}
		\State $V \gets Q Y$
		\Until{converged} \label{line:conv}
		\State \textbf{return} $\Call{SoughtEigenpairs}{V,\Lambda}$
		
		\EndFunction
	\end{algorithmic}
\end{algorithm}

Note, that in line \ref{line:reduced}, we solve a reduced (unbounded) eigenvalue problem of finding $N$ eigenvalues of $B=Q^H A Q \in \C^{N \times N}$, where $N$ denotes the upper bound on the number of eigenvalues of $A$ in $(-1,1)$. This can yield an enormously improved performance opposed to computing all the eigenvalues of $A \in \C^{n \times n}$, if $N$ is much smaller than $n$, the size of $A$.

There are two important notes to make on the implementation of this algorithm, for further understanding of Algorithm~\ref{algo:subspace_iteration}:
\begin{enumerate}
	\litem{Orthogonalisation} In line~\ref{line:ortho}, we compute a matrix $Q \in \C^{n \times N}$ with orthogonal columns, i.e.\ for two different columns $q_i, q_j$ of $Q$, we have $q_i^T q_j = 0$. To achieve this, the well-known \textit{\tool{QR} algorithm}\index{QR algorithm} can be used. A central feature of \FEAST{} is the absence of this orthogonalisation step, see \cite{feast,guettel}. This yields an improved performance, as the \tool{QR} algorithm is not efficient for large-scale orthogonalisation problems.
	
	\litem{Convergence criterion} For the convergence criterion in line~\ref{line:conv}, there are different approaches feasible. A common approach is monitoring the sum of absolute values of the computed eigenvalue approximations in $\Lambda=\text{diag}(\lambda_1,\dots, \lambda_k)$ lying in the search interval $(-1,1)$. This sum is known as \textit{eigentrace}\index{Eigentrace}. Thus, we have
	\begin{equation}
	\text{eigentrace}(\Lambda) = \sum_{\lambda_i \in (-1,1)} \abs{\lambda_i}.
	\end{equation}
	In each iteration, the algorithm computes the \textit{relative change} of the eigentrace compared to the previous iteration. The algorithm terminates when a user-selected tolerance is reached by the relative change. In the original publication of \FEAST{} \cite{feast}, a tolerance of $10^{-13}$ was proposed. More sophisticated approaches are possible, confer \cite{feast_conv}.
\end{enumerate}

\section{Convergence rate}
\label{sec:feast_conv}

In this thesis, we want to improve \FEAST{} by reducing its \textit{convergence rate}\index{FEAST!Convergence rate}\index{Convergence rate} for a variety of \HIEP{}s. In the following, we introduce this convergence rate. This result was shown in \cite{feast_subspace}.

Assume, the eigenvalues in the search interval $(-1,1)$ of $A$ were given by $\lambda_1, \dots, \lambda_k$, and all the eigenvalues of $A$ by $\lambda_1, \dots, \lambda_k, \dots, \lambda_n$, for some $k \in \N$. Then, the convergence rate of the \FEAST{} algorithm is
\begin{equation}
\label{eqn:conv}
\max_{1 \le i \le k}{\abbs{\frac{r(\lambda_{N+1})}{r(\lambda_i)}}},
\end{equation}
assuming the eigenvalues of $A$ can be ordered such that
\begin{equation}
\label{eqn:conv_required}
\abs{r(\lambda_1)} \ge \dots \ge \abs{r(\lambda_k)} \ge \dots \abs{r(\lambda_N)} \ge \abs{r(\lambda_{N+1})} \ge \dots,
\end{equation}
where $N$ denotes the upper bound on the eigenvalues of $A$ in $(-1,1)$ as in Algorithm~\ref{algo:subspace_iteration}.

The \FEAST{} algorithm converges faster, if this convergence rate is smaller. So, we want to reduce this convergence rate for various \textsc{Hermitian} matrices $A$ and eigenvalue distributions in this thesis.

The convergence rate explains, why we assumed \RFF{}s to be both real-valued and symmetric with respect to the $y$-axis in their definition. The former is a consequence of the eigenvalues of a \textsc{Hermitian} matrix being real. The latter is a consequence of an \RFF{} trying to reduce the convergence rate of \FEAST{}. If the underlying \CHIEP{} is not known, no assumption can made on the distribution of the eigenvalues. So, we assume, that eigenvalues not in the search interval lie equally-likely on either side of the search interval $(-1,1)$.

Note that assumption \eqref{eqn:conv_required} can always be assured by an appropriate choice of the \RFF{} $r$. If the \RFF{} has a larger absolute function value for values inside the search interval than outside\footnote{Any function with this property is called \textit{filter function}. We restrict out analysis to filter functions suitable for \FEAST{}, i.e.\ rational filter functions.}, assumption \eqref{eqn:conv_required} always holds. For instance, consider the function given by
\begin{equation}
\label{eqn:indicator}
r(x) = \begin{cases} 1, &\text{if } x \in (-1,1), \\ 
0, & \text{otherwise.} \end{cases}
\end{equation}
For this function, assumption \eqref{eqn:conv_required} is satisfied for any distribution of eigenvalues of $A$. More importantly, we note that this function induces the minimum convergence rate \eqref{eqn:conv}. It equals zero for any \textsc{Hermitian} matrix $A$.

Equation \eqref{eqn:indicator} is equivalent to the definition of the \textit{indicator function}\index{Indicator function} on $(-1,1)$, herein denoted as $\ind$. As we only examine the search interval $(-1,1)$, we can safely refer to this function as \textit{the} indicator function.

As the indicator function is discontinuous, it is impossible for a (continuous) \RFF{} to achieve this optimal convergence rate. So, we \textit{approximate}\index{Approximation} the indicator function by an \RFF{} in the following.

%% file: parts/RFF.tex
\preambletrans{A. Einstein}{Geometrie und Erfahrung}{Insofern sich die Sätze der Mathematik auf die Wirklichkeit beziehen, sind sie nicht sicher, und insofern sie sicher sind, beziehen sie sich nicht auf die Wirklichkeit.}{As far as theorems of mathematics do not refer to reality, they are not certain, and as far as they are certain, they do not refer to reality.}
\chapter{Existing rational filter functions}
\label{chp:rff}

\input{figures/RFF/gl}

In this chapter, we study three different approaches to approximating the indicator function by \RFF{}s. In other words, we analyse \RFF{}s for solving \CHIEP{}s using \FEAST{} efficiently. We discuss the two \RFF{}s primarily used by \FEAST{}, i.e. the \textsc{Gauss-Legendre} \RFF{} in Section~\ref{sec:gauss} and the \textsc{Zolotarev} \RFF{} in Section~\ref{sec:zolotarev}. In Section~\ref{chp:nonlin}, we introduce the \LSOP/ minimisation to further improve these \RFF{}s.

\section[Gauss-Legendre]{Using Gauss-Legendre quadrature}
\label{sec:gauss}

Our first approach to approximating the indicator function is exploiting the curve integral representation of the indicator function and solving the resulting integral numerically.\index{Numerical intergration} To achieve this, consider the unit circle $\gamma(t) := e^{\iu t}, \, t \in [0, 2 \pi]$. For ${x \in \R \setminus \set{-1,1}}$, the integral representation of the indicator function\footnote{For completeness, this is briefly proved in Lemma \ref{lem:indicator_integral}.} is given by
\begin{equation}
		\ind(x) = \left. \begin{cases} 1, &\text{if } x \in (-1,1), \\ 
		0, & \text{otherwise} \end{cases} \right\} = \frac{1}{2 \pi \iu} \, \int_{\gamma} \! \frac{1}{\zeta - x}  \, \mathrm{d}\zeta.
\end{equation}
Applying the definition of a curve integral, we obtain
\begin{equation}
	\ind(x) = \frac{1}{2 \pi \iu} \, \int_{0}^{2 \pi} \! \frac{\gamma'(t)}{\gamma(t) - x}  \, \mathrm{d}t =  \frac{1}{2 \pi} \, \int_{0}^{2 \pi} \! g_x(t) \, \mathrm{d}t,
\end{equation}
if we set
\begin{equation}
g_x(t) := \frac{e^{\iu t}}{e^{\iu t} - x}, \quad t \in [0, 2 \pi].
\end{equation}
Recall that $\Re$ denotes the real part of a complex number. Rewriting yields
\begin{equation}
	\label{eqn:integral}
	\ind(x) = \frac{1}{2 \pi} \, \Re \int_{0}^{\pi} \! g_x(t) + \conj{g_x(t)} \, \mathrm{d}t.
\end{equation}
As a simplification, we define the integrand of the resulting integral \eqref{eqn:integral} as
\begin{equation}
	h_x(t) := g_x(t) + \conj{g_x(t)}.
\end{equation}
With this notation, we compute the derived integral \eqref{eqn:integral} numerically. To achieve this, the so-called \textit{Gauss-Legendre quadrature} defines\index{Gauss-Legendre quadrature} \textit{nodes} $y_1, \dots, y_{2m} \in (0,\pi)$ to evaluate an integrand at and corresponding \textit{weights} $w_1, \dots, w_{2m} \in (0,\infty)$ to adding up the resulting values yielding
\begin{equation}
	\label{eqn:halfplane}
	\ind(x) = \frac{1}{2 \pi} \, \Re \int_{0}^{\pi} \! h_x(t) \, \mathrm{d}t  \approx \frac{1}{2 \pi} \Re \sum_{k=1}^{2m} w_k h_x(y_k).
\end{equation}
By symmetries of the nodes and weights in the \textsc{Gauss-Legendre} quadrature, the resulting rational function is an \RFF{} in terms of \eqref{eqn:rff}.
\pagebreak
\begin{definition}[Gauss-Legendre \RFF{}]
	\label{def:gauss}
	For $m\in \N$, we refer to the \RFF{} given by~\eqref{eqn:halfplane} as \textit{Gauss-Legendre} \RFF{}\index{RFF!Gauss-Legendre}\index{Gauss-Legendre RFF|see{RFF}} of degree $4m$.
\end{definition}
There exists no analytical representation of the weights and nodes of the \textsc{Gauss-Legendre} quadrature, only a numerical one \cite{dahmenreusken}. Also, it is possible to define the \textsc{Gauss-Legendre} \RFF{} on different paths rather than the unit circle. Still, the unit circle has proved to yield good results in practice. For a more extensive analysis see \cite{guettel,eigencount}.

A graphical illustration is given in Figure~\ref{fig:gaussianrff}. Note that the \textsc{Gauss-Legendre} \RFF{} does not have a sharp drop about $\pm 1$, but offers a good approximation of the indicator function everywhere else. This is a problem, if the non-sought eigenvalues are close the endpoints of the search interval, as of the convergence rate of \FEAST{} \eqref{eqn:conv}. In the following, we analyse an \RFF{} that attempts to solve this issue.

\section[Zolotarev]{Zolotarev, worst-case optimality}
\label{sec:zolotarev}

In this section, we present an approach to minimising the maximum absolute distance between the indicator function and a rational function. This can be phrased as the \textit{minimisation problem}\index{Minimisation problem}
\begin{equation}
	\label{eqn:minzolotarev}
	\argmin{r\in F}{\max_{x \in \R \setminus M} |\ind(x) - r(x)|},
\end{equation}
for a $G \in (0,1)$ and $M := [-G^{-1},-G] \cup [G,G^{-1}]$. The set $F$ denotes the rational functions of type $(2m,2m)$, for a fixed $m \in \N$.

Note that we choose some \textit{gap parameter}\index{Gap!Parameter} $G \in (0,1)$ to exclude a small interval enclosing $\pm 1$ in order to find a solution to the minimisation problem \eqref{eqn:minzolotarev}. This is necessary, because the indicator function is discontinuous at $\pm 1$ in contrast to \RFF{}s.

Recall the convergence rate of \FEAST{} as of \eqref{eqn:conv}. Intuitively, minimising the absolute difference between the indicator function and a rational function yields the optimal worst-case convergence rate, if we assume that there were no eigenvalues in the \textit{gap}s\index{Gap} as of $M = [-G^{-1},-G] \cup [G,G^{-1}]$.

The solutions to \eqref{eqn:minzolotarev} exist, can be computed numerically and yield an \RFF{}, as shown in \cite{guettel}, using results of \textsc{Zolotarev} from 1877.
\pagebreak
\begin{definition}[\textsc{Zolotarev} \RFF{}, Elliptic \RFF{}]
	\label{def:zolotarev}
	Let $m = 2n$ for an $n \in \N$. Then, we call the \RFF{} of degree $2m$ obtained from the version used in \FEAST{} that solves \eqref{eqn:minzolotarev} as \textit{Zolotarev \RFF{}}\index{RFF!Zolotarev}\index{Zolotarev RFF|see{RFF}} of degree $2m$. This \RFF{} is also known as \textit{Elliptic \RFF{}}.
\end{definition}

\input{figures/RFF/zolo}

The \textsc{Zolotarev} \RFF{} is illustrated in Figure~\ref{fig:ellipcticrff}. In contrast to the previously discussed \textsc{Gauss-Legendre} \RFF{}, the \textsc{Zolotarev} \RFF{} attains the same bound repeatedly. This is also depicted in Figure~\ref{fig:filters_log}. Another difference is that the decline about $\pm 1$ is sharper. Still, further away from the endpoints $\pm 1$ of the search interval $(-1,1)$, the \textsc{Gauss-Legendre} \RFF{} offers a better approximation of the indicator function.

\section[SLiSe]{SLiSe minimisation problems}
\label{chp:nonlin}\index{SLiSe}

In the two previous sections, we have introduced two approaches to gathering \RFF{}s: First, we rewrote the indicator function as an integral and derived a numerical solution via \textsc{Gauss-Legendre} quadrature. Second, we presented the \textsc{Zolotarev} \RFF{}, having a bounded absolute function value on the entire real line except a minor gap. In practice, these approaches are not flexible enough and do not offer sufficient degrees of freedom to suit different \HIEP{}s. To solve this issue, \textsc{Winkelmann} and \textsc{Di Napoli} introduced the concept of a \textit{weighted, squared approximation} of the indicator function \cite{jan}. They refer to their concept as \LSOP/. To illustrate this approach of \enquote{weighted, squared approximation}, we introduce the terms weighted, squared and approximation in the following.
\begin{enumerate}
	\litem{Implementing weighting} In particular, the function values of an \RFF{} about $\pm 1$ are interesting, when approximating the indicator function. This is, where the indicator function is discontinuous. So, we would like to put more weight on values about $\pm 1$, when approximating the indicator function. To formalise the term \textit{weight}, we introduce \textit{weight functions}\index{Weight function}. These functions shall yield non-negative values. The larger the weight function value, the more important is this point of the real line. We will require a weight function to be symmetric with respect to the $y$-axis. Thus, we ensure the same weighting for both positive and negative values.
	\vspace{5pt}
	\begin{definition}[Weight function]
		A function $\mathfrak{G}: \R \rightarrow [0,\infty)$ is called \textit{weight function}, if $\mathfrak{G}$ is symmetric with respect to the $y$-axis.
	\end{definition}
	A variety of weight functions were studied by \textsc{Winkelmann} and \textsc{Di Napoli} in \cite{jan}. Some prominent examples are given in Appendix~\ref{chp:weight_functions} for completeness.
	\litem{Measuring the weighted, squared difference} Recall, that in the previous chapter, we have introduced, how \RFF{}s are parametrised as per \eqref{eqn:rff}. Via integration, a weight function $\mathfrak{G}$ can be used to measure the approximation of the indicator function by an \RFF{} $\ratcp{r}{\beta,w}$ as per\index{Loss function}
	\begin{equation}
	\label{eqn:complexlossfunction}
	f(\beta, w) := \int_{-\infty}^{\infty} \! \mathfrak{G} (x) \, (\ind(x) - \ratcp{r}{\beta,w} (x))^2 \mathrm{d}x,
	\end{equation}
	for a fixed $m\in\N$ and all $\beta \in \C^m$, $w \in (\C \setminus \R)^m$.
	
	Note that $f$ is non-negative. Function values close to $0$ are supposed to indicate \enquote{well-performing} \RFF{}s, in contrast to function values further away from $0$.
	
	The used integral measures the overall \textit{squared} distance between a given \RFF{} and the indicator function. This is fundamentally different from the concept of the \textsc{Zolotarev} \RFF{} \eqref{eqn:minzolotarev}: Not the maximum absolute error of approximation is relevant, but the overall error.
	
	In practice, we can assume, that the integral in \eqref{eqn:complexlossfunction} always exists, because it is sufficient to consider weight functions $\mathfrak{G}$, that fulfil
	\begin{equation}
	\label{eqn:weight_radius}
	\mathfrak{G}(x) = 0, \quad \text{if} \ |x| \ge C,
	\end{equation}
	for some constant $C > 0$.
	
	We can make this assumption for two reasons: First, as a \textsc{Hermitian}\index{Hermitian} matrix $A$ has finitely many eigenvalues, there always exists an absolutely largest eigenvalue $\lambda$. This means that values absolutely larger than $\lambda$ do not require any weight, as of the convergence rate of \FEAST{} \eqref{eqn:conv}. We could choose a weight function, that fulfils \eqref{eqn:weight_radius} for $C:=\lambda +1$ to assure the existence of the studied integral in \eqref{eqn:complexlossfunction}. Second, any rational filter function $x \mapsto r(x)$ converges to $0$ for $x \rightarrow \pm \infty$ by \eqref{eqn:rff}.
	
	If the weight function $\mathfrak{G}$ is a \textit{step function}\index{Step function}, thus piece-wise constant, the function $f$ is twice continuously differentiable, as shown by \textsc{Winkelmann} and \textsc{Di Napoli} \cite{jan}. Moreover, in this case, the function $f$ and its gradient can be computed efficiently by matrix operations and solving standard definite integrals only. For the exact computational steps, see Theorem~\ref{thm:complexcostfunction}.
	
	\litem{Approximation by minimising the weighted, squared difference} For a given weight function $\mathfrak{G}$, the function $f$ measures the quality of overall approximation of the indicator function by an \RFF{}. So, we can use $f$ to approximate the indicator function by a \enquote{well-performing} \RFF{}. This can be phrased as the \LSOP/ minimisation problem\index{Minimisation problem}\index{SLiSe!RFF|see{RFF}}
	\begin{equation}
	\label{eqn:complex_minimiser}
	\argmin{\beta \in \C^m, w \in (\C \setminus \R)^m}{f(\beta, w)}.
	\end{equation}
	In contrast to the minimisation problem yielding the \textsc{Zolotarev} \RFF{}s \eqref{eqn:minzolotarev}, we do not measure the absolute difference, but rather the weighted \textit{squared} difference. This is a least-squares approximation, which comes close to the intuition of the average-case of the convergence rate of \FEAST{} \eqref{eqn:conv_required}: We do not minimise the absolute difference but rather the overall difference.
\end{enumerate}

As per \eqref{eqn:complexlossfunction}, the previous Equation \eqref{eqn:complex_minimiser} yields a \textit{non-linear} minimisation problem. Non-linear means that the problem cannot be described in terms of a linear equation system. This class of problems is computationally difficult. For linear problems, there exist such solvers such as the well-known \textit{\tool{QR} algorithm}. For non-linear problems, iterative solvers, such as \BFGS{}, can be employed. The performance of solving a non-linear minimisation problem depends on both the used solver and the minimisation problem itself. So, it would be desirable to gather access to a broad selection of non-linear minimisation algorithms. We accomplish this integration in the following chapter.

Integration into standard minimisation algorithms results in another benefit: Not only solutions of the minimisation problem \eqref{eqn:complex_minimiser} can be obtained more efficiently, but also modifications can be studied more rigorously for further improvement. We focus on such improvements to \LSOP/ in Chapter~\ref{chp:constraints}.

\pagebreak
Note that neither the function value $f$ \eqref{eqn:complexlossfunction} nor a solution to the minimisation problem \eqref{eqn:complex_minimiser} are connected with a small convergence rate of \FEAST{} as per \eqref{eqn:conv}. The quality of a solution of the minimisation problem \eqref{eqn:complex_minimiser} solely depends on the quality of the chosen weight function. However, there has not been discovered a criterion for the quality of weight functions in \FEAST{} in the past. They have to be obtained empirically. If a solution of the minimisation problem \eqref{eqn:complex_minimiser} could be obtained efficiently, a broad selection of weight function could be tested. With our well-performing integration into standard minimisation algorithms, we present how to obtain good weight functions algorithmically in Chapter~\ref{chp:tests}.

%% file: figures/RFF/gl.tex
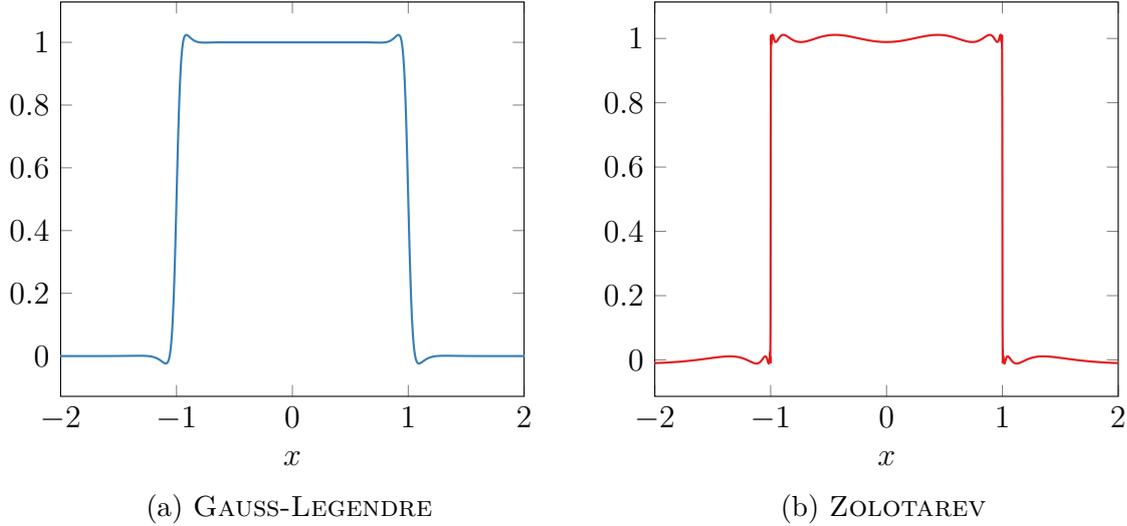
\begin{figure}
	\centering
	\begin{subfigure}[c]{0.48\textwidth}
		\tikzsetnextfilename{gaussianrff}
		\begin{tikzpicture}
		\begin{axis}[
		xmin=-2,
		xmax=2,
		xlabel={$x$}, 
		width=1\textwidth,
		height=0.3\textheight]
		\addplot[color=GAUSSC, thick] table[x index={0},y index={1}] {figures/RFF/rffs.dat};
		\end{axis}
		\end{tikzpicture}
		\subcaption{\textsc{Gauss-Legendre}}
		\label{fig:gaussianrff} 
	\end{subfigure}
	\begin{subfigure}[c]{0.48\textwidth}
		\tikzsetnextfilename{ellipcticrff}
		\begin{tikzpicture}
		\begin{axis}[
		xmin=-2,
		xmax=2,
		xlabel={$x$}, 
		width=1\textwidth,
		height=0.3\textheight]
		\addplot[color=ZOLOC, thick] table[x index={0},y index={1}] {figures/RFF/zolo.dat};		
		\end{axis}
		\end{tikzpicture}
		\subcaption{\textsc{Zolotarev}}
		\label{fig:ellipcticrff}
	\end{subfigure}
	\caption[\RFF{}s of \FEAST{}]{The discussed $16$-pole \RFF{}s used in \FEAST{}.}
	\label{fig:filters}
\end{figure}

%% file: figures/RFF/zolo.tex
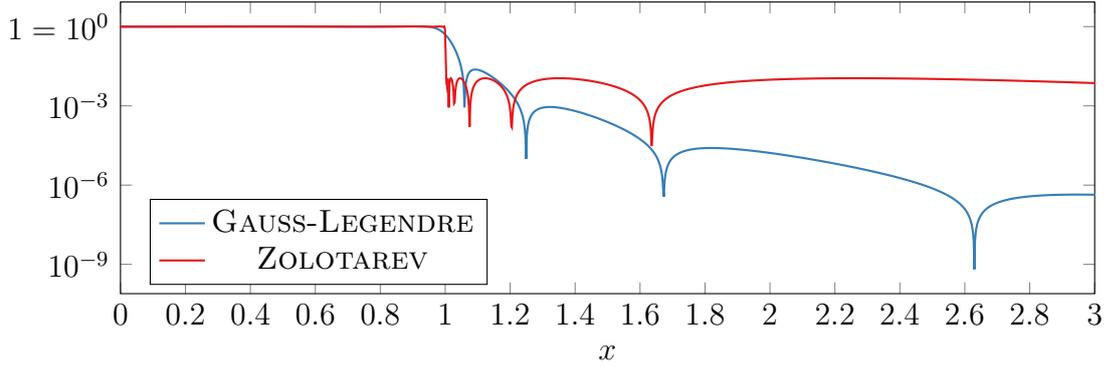
\begin{figure}
	\centering			
	\tikzsetnextfilename{filters_log}
	\begin{tikzpicture}
	\begin{axis}[
	xmin=0,
	xmax=3,
	ymode=log,
	ytick={0.000000001,0.000001,0.001,1},
	yticklabels={$10^{-9}$,$10^{-6}$,$10^{-3}$,$1 = 10^{0}$},
	legend pos=south west,
	xlabel={$x$}, 
	width=0.9\textwidth,
	height=0.24\textheight,
	]
	\addplot[color=GAUSSC, thick] table[x index={0},y expr=abs(\thisrowno{1})] {figures/RFF/rffs.dat};	
	\addplot[color=ZOLOC, thick] table[x index={0},y expr=abs(\thisrowno{1})] {figures/RFF/zolo.dat};	
	
	\addlegendentry{\textsc{Gauss-Legendre}}
	\addlegendentry{\textsc{Zolotarev}}
	\end{axis}
	\end{tikzpicture}
	\caption[Absolute value]{Logarithmic plot of the absolute value of the two discussed \RFF{}s, used in \FEAST{}, for the case of $16$ poles.}
	\label{fig:filters_log}
\end{figure}

%% file: parts/lossfunction.tex
\preamble{D. E. Knuth}{1974 Turing Award Lecture}{The real problem is that programmers have spent far too much time worrying about efficiency in the wrong places and at the wrong times; premature optimization is the root of all evil (or at least most of it) in programming.}
\chapter[Embedding into standard minimisation algorithms]{Using standard minimisation algorithms}
\label{sec:real}

\input{figures/lossfunction/quadratic}

In the previous chapter, we have introduced a non-negative, real-valued, twice continuously differentiable function, that measures the overall approximation of the indicator function by an \RFF{}, denoted
\begin{equation*}
f(\beta, w), \quad \text{for all} \ \beta \in \C^m, w \in (\C \setminus \R)^m,
\end{equation*}
as per \eqref{eqn:complexlossfunction}, for a fixed $m \in \N$.

This motivated the minimisation problem \eqref{eqn:complex_minimiser} of $f$ to gather \enquote{well-performing} \RFF{}s for \FEAST{}. We now want to solve this minimisation problem\index{Minimisation problem} numerically, using standard minimisation algorithms, such as \BFGS{}, to obtain well-performing \RFF{}s efficiently. For this, we study the central input of every minimisation algorithm: The function to minimise.

\section{Loss function}

Most of the standard minimisation algorithms necessitate a standard input format for the function to minimise: the \textit{loss function}. This is defined as follows.
\begin{definition}[Loss function, residual]
	\label{def:loss}\index{Loss function}\index{Residual}
	Let $n\in\N$. Then, a twice continuously differentiable function $g: \R^n \rightarrow [0,\infty)$ is called \textit{loss function} or \textit{loss function of $n$ variables}. A function value of a loss function is called \textit{residual}.
\end{definition}
The term \enquote{residual} is derived from the application of loss functions in minimisation problems: If we have found a solution to a minimisation problem induced by a loss function (thus a minimiser of the loss function), the function value of the solution (the function value of the minimiser) is the remainder, the \textit{residual}, in the minimisation.
\begin{example}
	\label{exmpl:loss_function}
	The previous Definition~\ref{def:loss} generalises the \textit{quadratic loss function}\index{Loss function!Quadratic}, i.e.\ a function $\R^n \rightarrow \R, \ x \mapsto x^T A x$, for some fixed matrix $A \in \R^{n \times n}$ and $n \in \N$. If A is \textit{positive-semidefinite}, this is equivalent to $A$ being non-negative. In this case, we have a loss function in our terms. Two examples are depicted in Figure~\ref{fig:quadratic}. Note that a loss function does not necessarily have a unique minimiser. This is the case for non-convex quadratic loss functions.
\end{example}
The example of quadratic loss function is to illustrate loss functions. In our case, we have to handle by far more complex examples of loss functions.

\section{Transformation to real arguments}

Our function $f$ almost matches the Definition~\ref{def:loss} of a loss function. The only difference is that our function requires complex arguments. We perform a transformation to real arguments in this section. To illustrate our approach, consider the following example first.

\begin{example}
	Let $g$ be the real-valued function with complex arguments, given by
	\begin{equation}
	g(z) := \conj{z} z = \Re(z)^2 + \Im(z)^2, \quad \text{for} \ z \in \C.
	\end{equation}
	Its gradient is $\nabla g(z) = \conj{z}$. This can be seen by \enquote{treating} $\conj{z}$ as a constant and differentiating with respect to $z$.
	
	Any $z\in\C$, can be written as $z = x + \iu y$ for $x,y \in \R$. So, we have $g(z) = g(x + \iu y) = x^2 + y^2$ for $x,y \in \R$. This is a function mapping reals to reals. We can equivalently define
	\begin{equation}
	\tilde{g}(\colvecalt{x}{y}) := x^2 + y^2 = g(x + \iu y).
	\end{equation}
	This is a loss function as per Definition~\ref{def:loss}. We have introduced this function before, confer Figure~\ref{fig:quadratic1} and the previous Example~\ref{exmpl:loss_function}. The gradient of $\tilde{g}$ is 
	\begin{equation}
	\nabla \tilde{g}(\colvecalt{x}{y}) = 2x + 2y = 2 \conj{\nabla g (x + \iu y)}.
	\end{equation}
	We have pointed out a connection between $g$ and $\tilde{g}$. This result is important, as it will be generalised in the following.
\end{example}

For the general case, it was shown in \cite{complexoptimization}, that we can transform our loss function to real arguments by \enquote{splitting} the arguments into the real and imaginary parts. A function $g:\C^n \rightarrow \R$ can be transformed to an $\tilde{g}:\R^{2n} \rightarrow \R$ as per
\begin{equation}
\label{eqn:complex_transform}
\tilde{g}(\colvecalt{x}{y}) := g(x + \iu y),
\end{equation}
for $x,y\in \R^n$. According to \cite{complexoptimization}, the gradient is given by
\begin{equation}
\label{eqn:complex_transform2}
\nabla \tilde{g} (\colvecalt{x}{y}) = 2 \conj{\nabla g (x + \iu y)},
\end{equation}
for $x,y\in \R^n$. We can now safely call our function $f$ loss function\index{SLiSe!Loss function}, because it can be transformed to such, as per this procedure. The function values of $f$ will be called \textit{residuals}\index{SLiSe!Residual}\index{Residual} as in the case of loss functions. Our loss function can now be minimised using standard minimisation algorithms. In the following, we present an example of a powerful standard minimisation algorithm.
\section{BFGS minimisation algorithm}
\label{sec:bfgs}

As a minimisation algorithm, we propose the \textit{Broyden-Fletcher-Goldfarb-Shanno} (\BFGS{})\index{BFGS} algorithm \cite{numericaloptmization}. We discuss the algorithm in the following, while comparing to other prominent minimisation algorithms.\index{Minimisation problem}

First of all, \BFGS{} requires a loss function $g$ as per Definition~\ref{def:loss}, its gradient and a \textit{starting point}\index{Starting point}\index{BFGS!Starting point}\index{RFF!Starting} $x_0$ as inputs. For $k=1,2,\dots$, the algorithm approximates a minimiser of $g$ by \textit{iterates}\index{Iterate}\index{BFGS!Iterate} $x_k \in \R^n$, such that
\begin{equation}
g(x_0) > g(x_1) > \dots > g(x_k) \ge 0.
\end{equation}
This is supposed to converge to a minimiser of $g$. Note that $g(x) \ge 0$ for $x \in \R$ by the Definition~\ref{def:loss} of a loss function.

The \BFGS{} algorithm has multiple properties that make it an appropriate choice compared to other minimisation algorithms in our case. These properties are discussed in the following.
\begin{itemize}
	\litem{Efficiency} The algorithm does not require any linear system solves. This makes it superior amongst other minimisation algorithms, such as \textit{Levenberg-Marquardt}\index{Levenberg-Marquardt}. Furthermore, the algorithm memorises the history of iterates, in contrast to other methods, such as \textit{Gradient descent}\index{Gradient descent}\index{Projected!Gradient descent}. Tracking the iterates enables \BFGS{} to converge to a minimiser efficiently.
	
	\litem{Improvement} The algorithm requires a starting point as an input. Because we have already gathered prominent examples of \RFF{}s in Chapter~\ref{chp:rff}, we can use those as a starting point. This is supposed to yield an improved \RFF{} compared to the original \RFF{} \cite{jan}.
	
	\litem{Extensibility} By further conditions on the minimisation problem, the possible solutions can be constrained to obtain different \RFF{}s. The \BFGS{} algorithm is a versatile tool for embedding such conditions. In the following chapter, we study such conditions.
\end{itemize}

\paragraph{Extensive analysis of BFGS.} We provide an extensive analysis of the \BFGS{} algorithm in Appendix~\ref{chp:bfgs}. \BFGS{} is not a trivial algorithm, as it incorporates a multitude of advances in minimisation research \cite{armijo,hagerzhang,moerthuente,numericaloptmization,wolfe1,wolfe2}.

\section{Performance of BFGS}
\label{sec:bfgs_performance}

To demonstrate the performance of our approach, we compare \BFGS{} to the approach in the original publication by \textsc{Winkelmann} and \textsc{Di Napoli} in the case of \textit{box-constrained}\index{Box-constraint} \LSOP/ \cite{jan}. We have not introduced box-constrained \LSOP/ yet. This is the aim of the following chapter. However, we refer to the same set-up as in the original paper, which justifies discussing the performance of \BFGS{} in this case upfront.

The \textit{projected gradient descent}\index{Projected gradient descent} algorithm was used by \textsc{Winkelmann} and \textsc{Di Napoli} \cite{jan}. We use the \LBFGSB{}\index{BFGS!L-BFGS-B}\index{L-BFGS-B|see{BFGS}} algorithm, an extension of \BFGS{} \cite{lbfgs1}. As in the case of \BFGS{}, projected gradient descent does not involve linear system solves. In contrast to \BFGS{}, projected gradient descent does not memorise the history of iterates. This results in worse performance. For a more extensive analysis of projected gradient descent in the context of \RFF{}s, see the original publication \cite{jan}. For the \LBFGSB{} algorithm, we used the \tool{FORTRAN} implementation by the authors of the algorithm \cite{lbfgs1, lbfgs2,lbfgs3}. We employed the most recent version \tool{3.0} in our experiments.\footnote{We complied \LBFGSB{} using the supplied \tool{MAKEFILE} that builds upon \tool{gcc}, used in version \tool{5.4.0}. The processor was an \tool{Intel Core i5-4200M}. In \LBFGSB{}, we used standard parameters for maximum accuracy as per the documentation \cite{lbfgs2}, given by \tool{factr = 1e1}, \tool{m = 20}.} For projected gradient descent, we used the implementation of \textsc{Winkelmann} and \textsc{Di Napoli} who kindly granted access.

For a comparison, we refer to the same set-up as in the original publication for box-constraints by \textsc{Winkelmann} and \textsc{Di Napoli} first. Especially, this means that we employed the weight function Box-\LSOP/ (given in Example~\ref{exmpl:box_slise}) and the $16$-pole \textsc{Zolotarev} \RFF{}\index{RFF!Zolotarev} as the starting point of the minimisation algorithm. This number of poles is the standard choice in \FEAST{}.\index{FEAST}

In the case of projected gradient descent, the resulting residual\index{Residual} is $7.57\cdot 10^{-4}$. In the case of \LBFGSB{}, the resulting residual is $4.72\cdot 10^{-4}$, which is about $38$\% smaller. Thus, our approach yields a better \RFF{} in terms of weighted, squared distance. Both algorithms require less than $400$ function evaluations, which is a good performance for practical use.

This raises the question, why the algorithms converge to different residuals and \RFF{}s. As \LSOP/\index{SLiSe} is a non-convex problem, there may exist different minimisers having different residuals \cite{jan}. However, this is not the case, as the resulting \RFF{} of projected gradient descent can be further minimised to a residual of $4.72\cdot 10^{-4}$ by \LBFGSB{}. This is the same residual as in the case of \LBFGSB{} before. In fact, projected gradient descent does only slightly minimise the starting residual. Whereas the resulting residual of projected gradient descent is $7.57\cdot 10^{-4}$, the residual of the starting \textsc{Zolotarev} \RFF{} is $8.09\cdot 10^{-4}$. So, projected gradient descent does not offer the same accuracy as \LBFGSB{}.

\input{figures/experiments/box_residuals}

Also, the case is different, when another starting point is used. For the $16$-pole \textsc{Gauss-Legendre} \RFF{}\index{RFF!Gauss-Legendre} as the starting point, the results are depicted in Figure~\ref{fig:residuals}. On the $x$-axis, the number of evaluations of the loss function  $f$ after each iteration of the minimisation algorithm is denoted. On the $y$-axis, the residual after each iteration of the algorithm is denoted.

We find that the \LBFGSB{} algorithm needs about $500$ loss function evaluations, compared to about $900,000$ in the case of projected gradient descent. Furthermore, despite fewer loss function evaluations, \LBFGSB{} yields a smaller residual of $4.72\cdot 10^{-4}$ opposed to $7.08\cdot 10^{-4}$ for projected gradient descent. \LBFGSB{} results in the same residual and only about $100$ function evaluations more as for the \textsc{Zolotarev} \RFF{} as the starting point. In contrast, projected gradient descent results in a better (smaller) residual compared to before, but a much larger number of function evaluations. Nevertheless, the resulting residual of projected gradient descent is not better (smaller) than in the case of \LBFGSB{}.

We have also tested different number of poles and different weight functions. However, we found similar results in these cases. We conclude that our approach significantly increases the resilience of minimisation against the choice of parameters in the minimisation process. \LBFGSB{} outperforms projected gradient descent in both performance and accuracy.

%% file: figures/lossfunction/quadratic.tex
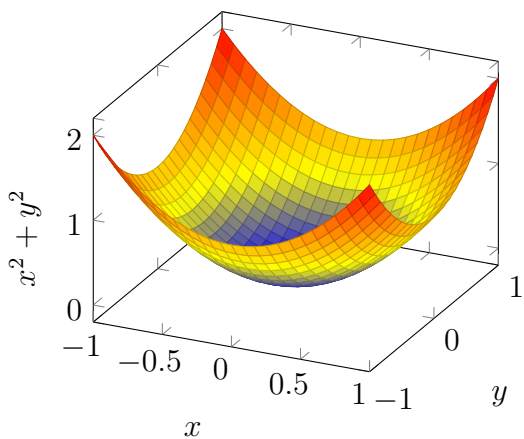
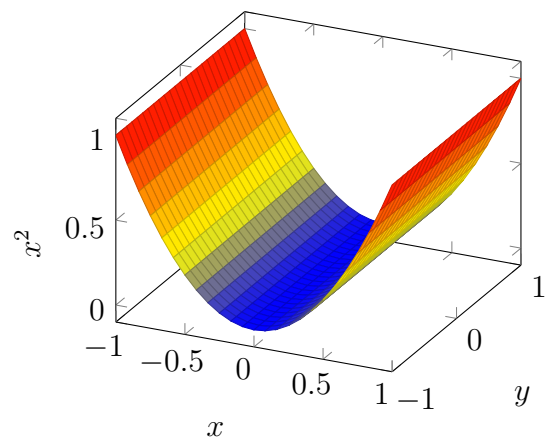
\begin{figure}[H]
	\centering
	
	\begin{subfigure}[c]{0.48\textwidth}
		\tikzsetnextfilename{quadratic1}
		\begin{tikzpicture}
		\begin{axis}[
		domain=-1:1,y domain=-1:1,
		xlabel={$x$}, ylabel={$y$}, zlabel={$x^2 + y^2$},
		width=.9\textwidth,
		height=0.28\textheight]
		\addplot3[surf] {x^2 + y^2};
		\end{axis}
		\end{tikzpicture}
		\subcaption{Convex quadratic loss function.}
		\label{fig:quadratic1}
	\end{subfigure}
	\begin{subfigure}[c]{0.48\textwidth}
		\tikzsetnextfilename{quadratic2}
		\begin{tikzpicture}
		\begin{axis}[
		domain=-1:1,y domain=-1:1,
		xlabel={$x$}, ylabel={$y$}, zlabel={$x^2 $},
		width=.9\textwidth,
		height=0.28\textheight]
		\addplot3[surf] {x^2};
		\end{axis}
		\end{tikzpicture}
		\subcaption{Non-convex quadratic loss function.}
		\label{fig:quadratic2}
	\end{subfigure}
	\caption[Loss functions]{Examples of quadratic loss functions.}
	\label{fig:quadratic}
\end{figure}

%% file: figures/experiments/box_residuals.tex
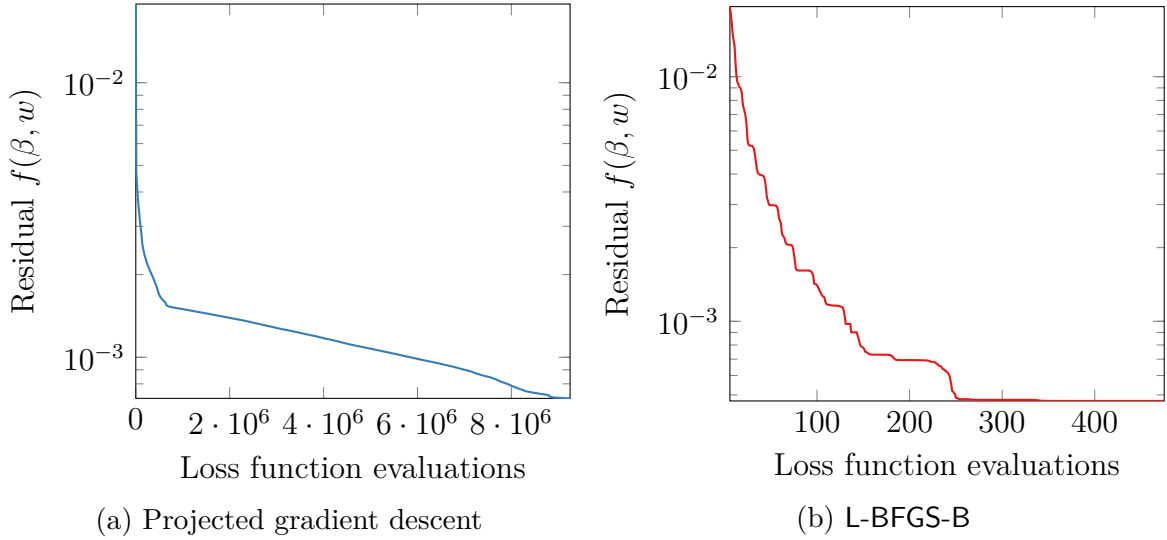
\begin{figure}
	\centering
	\begin{subfigure}[c]{0.48\textwidth}
		\tikzsetnextfilename{projected_descent_residuals}
		\begin{tikzpicture}
			\begin{axis}[
			scaled x ticks=false, 
			xlabel={Loss function evaluations}, ylabel={Residual $f(\beta,w)$},
			enlargelimits=false,
			ymode=log,
			width=0.95\textwidth,
			height=0.3\textheight,
			]
			\addplot[color=GAUSSC, thick] table[x index={2},y index={1}] {figures/experiments/projected_descent_residuals.dat};
			\end{axis}
		\end{tikzpicture}
		\subcaption[Projected gradient descent]{Projected gradient descent}
		\label{fig:projected_descent_residuals}
	\end{subfigure}
	\begin{subfigure}[c]{0.48\textwidth}
		\tikzsetnextfilename{lbfgsb_residuals}
		\begin{tikzpicture}
		\begin{axis}[
		xlabel={Loss function evaluations}, ylabel={Residual $f(\beta,w)$},
		enlargelimits=false,
		ymode=log,
		width=0.95\textwidth,
		height=0.3\textheight,
		]
		
		\addplot[color=ZOLOC, thick] table[x index={2},y index={1}] {figures/experiments/lbfgsb_residuals.dat};
		\end{axis}
		\end{tikzpicture}
		\subcaption{\LBFGSB{}}
		\label{fig:lbfgsb_residuals}
	\end{subfigure}
	\caption[Residuals in box-constraints]{Comparison of \LBFGSB{} and projected gradient descent in box-constrained \LSOP/.}
	\label{fig:residuals}
\end{figure}

%% file: parts/constraints.tex
\preamble{C. A. R. Hoare}{1980 Turing Award Lecture}{There are two ways of constructing a software design: One way is to make it so simple that there are obviously no deficiencies and the other way is to make it so complicated that there are no obvious deficiencies.}
\chapter{Constrained rational filter functions}
\label{chp:constraints}
In the preceding chapter, we have derived how to embed \LSOP/ into standard minimisation algorithms. As a prominent example of such, we studied the \BFGS{} algorithm. This opens the opportunity to make use of further features of such minimisation algorithms: Adding \textit{constraints}\index{Constraint}. This means that we are not interested in \textit{any} \RFF{} in \LSOP/ minimisation, but rather in one subject to certain constraints to further improve \RFF{}s. In other words, we study modifications to the original \LSOP/ minimisation problems \eqref{eqn:complex_minimiser}. The \BFGS{} algorithm allows to implement constraints relatively easily.

In the following, we introduce three different types of constraints: In Section~\ref{sec:box}, we improve the performance of performing a single \FEAST{} iteration, when \LSOP/ \RFF{}s are used. In Section~\ref{sec:shape}, we present constraints to shape \RFF{}s. In Section~\ref{sec:prob}, we constrain \RFF{}s using a given probability distribution to improve the performance of \RFF{} for a particular set of \HIEP{}s.

\section{Box-constraints}
\label{sec:box}
\subsection{Problem}
Recall that the usage of \RFF{}s in \FEAST{} results in solving linear systems of the form
\begin{equation}
\label{eqn:feast_rff2}
\alpha (A - z I)^{-1} v,
\end{equation}
for vectors $v \in \C^n \setminus \set{0}$ and $\alpha \in \C$, $z \in \C \setminus \R$, as per \eqref{eqn:feast_rff}.

This is a problem, if a pole $z$ of an \RFF{} approaches an eigenvalue of $A$ in the minimisation process. In this case, the \textit{condition number}\index{Condition number} of the matrix $\alpha (A - z I)$ becomes very large. This yields slow convergence for iterative linear system solvers in \FEAST{}.\index{FEAST}

In the following, we analyse this condition number further. We neglect the scalar $\alpha$, as it is a constant factor in this analysis. The condition number of the remainder is given by\footnote{Herein, we analyse the condition number with respect to the $2$-norm. This is a natural choice in the context of complex numbers. Moreover, the matrix $A$ is assumed to be \textsc{Hermitian}, which implies that $A - z I$ is normal. In this case, our claim holds.}
\begin{equation}
\text{cond}(A - z I) = \frac{\max_i\abs{\lambda_i-z}}{\min_j\abs{\lambda_j-z}},
\end{equation}
where $\lambda_1, \dots, \lambda_n$ denote the eigenvalues of $A$.

We rewrite and estimate 
\begin{equation}
\text{cond}(A - z I) = \frac{\max_i\abs{\lambda_i-z}}{\min_j\abs{\lambda_j-\Re(z) - \iu \Im(z)}} \le \frac{\max_i\abs{\lambda_i-z}}{\abs{\Im(z)}},
\end{equation}
as all eigenvalues $\lambda_j$ of the \textsc{Hermitian} matrix $A$ are real. Recall that $\Re$ denotes the real, respectively $\Im$ the imaginary part.

We assume $\abs{z} < 1$, as we will only consider \RFF{}s with such poles $z$. This is related to the search interval being $(-1,1)$ in a \CHIEP{}. Then, we can further simplify using the triangle inequality
\begin{equation}
\text{cond}(A - z I) < \frac{\max_i |\lambda_i| + 1}{ \abs{\Im(z)}}.
\end{equation}
As $\max_i |\lambda_i| + 1$ is a constant factor, we neglect this scalar and define the following quantity.
\begin{definition}[Worst-case condition number]\index{FEAST!Worst-case condition number}
	\label{def:wccn}
	Let $r$ be an \RFF{}. Then, the \textit{worst-case condition number} in \FEAST{} is (up to a constant factor)
	\begin{equation}
	\text{cond}_\text{worst}(r) := \max_\text{$z$ is pole of $r$}{\frac{1}{ \abs{\Im(z)}}}.
	\end{equation}
\end{definition}
\subsection{Solution}
Fix an \RFF{} $\ratcp{r}{\beta,w}$ for some $\beta \in \C^m, w\in (\C \setminus \R)^m$ and $m \in \N$, parametrised as in \eqref{eqn:rff}. Recall that the values $w_1, \dots, w_m$ characterise the poles of the \RFF{} $\ratcp{r}{\beta,w}$ up to symmetries.

To limit the worst-case condition number, we introduce a \textit{lower bound}\index{Lower bound|see{Box-constraint}} $lb > 0$ on the absolute imaginary parts of the poles of the \RFF{} $\ratcp{r}{\beta,w}$ as per
\begin{equation}
\abs{\Im(w_i)} \ge lb, \quad \text{for} \ i=1,\dots,m.
\end{equation}
In this case, as of Definition~\ref{def:wccn}, we have
\begin{equation}
	\text{cond}_\text{worst}(\ratcp{r}{\beta,w}) \le \frac{1}{lb}.
\end{equation}
For a given lower bound $lb > 0$, the previously studied minimisation problem \eqref{eqn:complex_minimiser} can be reformulated as
\begin{subequations}
	\label{eqn:box}\index{Minimisation problem}
	\begin{equation}
	\argmin{\beta \in \C^m, w \in (\C \setminus \R)^m}{f(\beta, w)},
	\end{equation}
	\begin{equation}
	\text{subject to:} \ \abs{\Im(w_i)} \ge lb, \quad \text{for} \ i=1,\dots,m,
	\end{equation}
\end{subequations}
as introduced in \cite{jan}. In a general minimisation environment, such bounds are known as \textit{box-constraints}\index{Box-constraint}. This is demanding arbitrary lower and upper bounds on the components of the sought solution. For the solution of box-constrained problems, there exists a well-known extension to \BFGS{}, the \LBFGSB{}\index{BFGS!L-BFGS-B}\index{BFGS!Constrained} algorithm \cite{lbfgs1}.

For more background on box-constrained minimisation in the case of \LBFGSB{}, see Appendix~\ref{sec:constrained_optimisation}.
\clearpage
\subsection{Performance}
We have previously discussed the performance of our implementation of box-constraints in Section~\ref{sec:bfgs_performance}. Compared to the current implementation of box-constraints in \LSOP/, our \LBFGSB{}\index{BFGS!L-BFGS-B} approach has significantly improved both performance and accuracy. In this subsection, we discuss the impact of different lower bounds on both the residual and the worst-case condition number in terms of Definition~\ref{def:wccn} of \RFF{}s.

We employ the same set-up as in the original \LSOP/ publication \cite{jan} for box-constrained minimisation: We use the Box-\LSOP/ weight function in the minimisation problem. This weight function was explicitly designed for box-constrained minimisation. Moreover, we choose the $16$-pole \textsc{Zolotarev}\index{RFF!Zolotarev} \RFF{} as the starting point of the minimisation for two reasons. First, the number of $16$ poles results is the default value in \FEAST{}\index{FEAST}, as it performs well in practice \cite{feast}. Second, the minimum absolute imaginary of the poles of this \RFF{} is only about $0.0022$, which is low compared to other \RFF{}s. This can be seen as follows: Using the notation as of \eqref{eqn:rff}, i.e.\ only employing the poles in the upper-left complex half-plane, we can denote the $16$-pole \textsc{Zolotarev} \RFF{} $\ratcp{r}{\beta,w}$ as
\begin{center}
	\begin{tabular}{l|l}
		\hline
		\normalsize{Poles \(w\)} & \normalsize{Coefficients \(\beta\)}\\
		\hline
		\(-0.999998+\textbf{0.0021993} \iu\)& \(0.00089892-0.00000198 \iu\)\\
		\(-0.999851+0.0172345 \iu \)& \(0.00524579- 0.0000904\iu\)\\
		\(-0.993336+0.115256 \iu \)& \(0.0346254-0.00401754 \iu\)\\
		\(-0.739835+0.672789 \iu \)& \(0.150517-0.136877 \iu\)\\
		\hline
	\end{tabular}
\end{center}
For a comparison to the other studied \RFF{}s, see Appendix~\ref{chp:resulting_rffs}. There, these \RFF{}s are denoted.

As minimisation algorithm, we use \LBFGSB{}\index{BFGS!L-BFGS-B} in the same set-up as in Section~\ref{sec:bfgs_performance}. We have computed both the worst-case convergence rates, as well as the residuals for $100$ equidistant lower bounds in $[0,0.0022]$. In the case of the \textsc{Zolotarev} \RFF{}, lower bounds larger than $0.0022$ are not feasible, as \LBFGSB{} requires that the starting point lies in the bounds.

\input{figures/experiments/different_lower_bounds_residuals}

The results are depicted in Figure~\ref{fig:different_lower_bounds_residuals}. For $lb=0.0022$, we obtain a residual of about $4.72\cdot10^{-4}$\index{Residual} and a worst-case condition number of about $6.13 \cdot10^{5} $. Conversely, for $lb=0.0011$, we find a residual of about $3.80\cdot10^{-4}$ and a worst-case condition number of about $17.8 \cdot10^{5}$, which is thrice as much.

For lower bounds smaller than $0.0011$, neither the residuals nor the worst-case condition numbers change in a meaningful way. The minimum residual is already reached by the minimisation algorithm for this lower bound.

We note, that we have obtained almost continuous functions of residuals and worst-case convergence rates. This means, that the chosen lower bound has a direct impact on both these values, as long as the minimiser of the original, unbounded \LSOP/ minimisation problem has not been reached yet.

We conclude that limiting the lower bound can reduce the worst-case condition number significantly. In other words, limiting the lower bound can indeed improve the performance of iterative linear system solves within \FEAST{}. 

\section{Shape constraints}
\label{sec:shape}
We have seen examples of \RFF{}s before, i.e.\ the \textsc{Gauss-Legendre} and the \textsc{Zolotarev} \RFF{}. They share a feature, that is desirable in practice: They have \textit{non-increasing extrema}\index{RFF!Non-increasing extrema} outside the search interval $(-1,1)$. This matches the intuition of an \RFF{} quite well: An \RFF{} delivers fast convergence in \FEAST{}, if the eigenvalues outside the search interval are not close its endpoints. An \RFF{} with a different behaviour may contradict the intuition of the \FEAST{} user.

There exist \LSOP/ \RFF{}s having \textit{increasing} extrema. An example of such is depicted in Figure~\ref{fig:nondecreasing_log}. The \textsc{Gauss-Legendre} \RFF{} has non-increasing extrema, whereas the \enquote{Increasing \RFF{}} has a local extremum at about $2$ with a larger absolute value than the extrema inside $[0,2)$.

\input{figures/constraints/nondecreasing}

In the following, we present a method to prevent such \RFF{}s in \LSOP/ minimisation. To achieve this, we introduce a set of constraints to be added to the original minimisation problem \eqref{eqn:complex_minimiser} in the following.

Let $\beta_0 \in \C^m, w_0 \in (\C \setminus \R)^m$ characterise the starting \RFF{} $\ratcp{r}{\beta_0,w_0}$ used in the minimisation algorithm, e.g.\ \BFGS{} as per Section~\ref{chp:nonli}. Choose $k \in \N$ \textit{evalutation points}, ${1 < x_1 < \dots < x_k}$. Then for each evalutation point $x_i$, we add a constraint to the minimisation problem as follows
\begin{equation}
	\label{eqn:shape_constraint1}
	\abs{\ratcp{r}{\beta,w}(x_i)} \le (1+c) \cdot \abs{\ratcp{r}{\beta_0,w_0}(x_i)},
\end{equation}
where $c \in (0,1)$ controls the percentage of deviation of the resulting \RFF{} from the starting \RFF{} at $x_i$.

Note, that it is sufficient to choose a small number of evaluation points $x_i$ near the endpoints of the search interval: On the one hand, an \RFF{} is a continuous function, which limits the steepness of change in function value. However, the larger the degree of the sought \RFF{}, the more evaluation points and constraints are required, as such \RFF{}s offer more degrees of freedom. On the other hand, an \RFF{} converges to zero for both large positive and negative values as pointed out before.

\pagebreak
We note that the values on the right hand sides of the in-equations \eqref{eqn:shape_constraint1} are constant. So, we can rewrite
\begin{equation}
	\label{eqn:shape_constraint}
	\abs{\ratcp{r}{\beta,w}(x_i)} \le C_i,
\end{equation}
defining the constants
\begin{equation}
	C_i := (1+c)  \cdot \abs{\ratcp{r}{\beta_0,w_0}(x_i)},
\end{equation}
for a control parameter $c \in (0,1)$ and $i=1,\dots,k$.

Faster convergence in minimisation algorithms can be achieved, if both the loss function and the constraints are differentiable; \RFF{}s are, but the absolute value is not. Instead of \eqref{eqn:shape_constraint}, we could add the following two constraints to the minimisation problem
\begin{subequations}
	\begin{align}	
	\ratcp{r}{\beta,w}(x_i) &\le C_i, \\
	-\ratcp{r}{\beta,w}(x_i) &\le C_i,
	\end{align}
\end{subequations}
for $i=1,\dots,k$.

We will refer to this type of constraints as \textit{shape constraints}\index{Shape constraint}, because they enforce the shape of the sought \RFF{} in the minimisation process.

The resulting minimisation problem can be phrased as
\begin{subequations}
	\label{eqn:shape_constraint_minimisation}\index{Minimisation problem}
	\begin{equation}
	\argmin{\beta \in \C^m, w \in (\C \setminus \R)^m}{f(\beta, w)},
	\end{equation}
	\begin{equation}
	\text{subject to:} \ \pm \ratcp{r}{\beta,w}(x_i) \le C_i, \quad \text{for} \ i=1,\dots,k.
	\end{equation}
\end{subequations}
For instance, this minimisation problem can be solved using the \tool{NLOpt} minimisation package \cite{nlopt}. This package provides multiple suitable minimisation algorithms, e.g.\ the \BFGS{}-based\index{BFGS} \tool{SLSQP} algorithm \cite{SLSQP}. In the following, we will not analyse this type of constraints any further, as they do not improve convergence of \FEAST{}\index{FEAST}, as per the definition of its convergence rate \eqref{eqn:conv}. Moreover, these constraints may limit the convergence of the \FEAST{} algorithm, as every added constraint potentially leads to an \RFF{} having a larger residual in \LSOP/ minimisation. This means that a resulting \RFF{} would be supposed to perform worse in \FEAST{}. However, for gathering \RFF{}s ready for production and the end-user, this technique might be considered.
\section{Probabilistic rational filter functions}
\label{sec:prob}
The convergence rate of \FEAST{} \eqref{eqn:conv} determines the performance of \FEAST{}. Previously, we have assumed to have no prior knowledge of the underlying eigenvalue distribution, in the design of \RFF{}s to suit a variety of \CHIEP{}s.

\input{figures/experiments/exp_conv}

There are cases, where we have prior knowledge of the probability distribution of the eigenvalues, for instance in \textit{self-consistent field calculations} from material science \cite{feast}. In this case, a probability distribution can be employed to obtain \textit{Probabilistic \RFF{}s}\index{RFF!Probabilistic} as follows.

Assume a continuous probability distribution $D$ and its corresponding probability density function $h : \R \rightarrow [0,\infty)$ that indicates the eigenvalue distribution of the studied \textsc{Hermitian} matrix $A$. As in the case of the \textsc{Zolotarev} \RFF{} in Section \ref{sec:zolotarev}, we assume a \textit{gap parameter}\index{Gap!Parameter} $G \in (0,1)$.

For values \textit{inside} the search interval, we define
\begin{equation}
	I := [-G,G].
\end{equation}
For values \textit{outside} the search interval, we introduce
\begin{equation}
	O := [-\infty,-G^{-1}]\cup[G^{-1},\infty].
\end{equation}
For given independent random variables $X,Y$ distributed according to $D$ and an \RFF{} $\ratcp{r}{\beta,w}$, we can determine the \textit{expected} convergence rate as of
\begin{equation}
\label{eqn:expect}
E(\beta,w) := \expect*{\abs{\dfrac{\ratcp{r}{\beta,w}(Y)}{\ratcp{r}{\beta,w}(X)}}| X \in I, \, Y \in O},
\end{equation}
parametrising an \RFF{} by $\beta \in \C^m, w \in (\C \setminus \R)^m$ and $m \in \N$ as per \eqref{eqn:rff}.

For independent random variables, the \textit{conditional expectation}\index{Conditional Expectation} \eqref{eqn:expect} is defined as
\begin{subequations}
	\begin{align}	
	E(\beta,w) =& \frac{1}{P[X \in I] \cdot P[Y\in O]} \int_{I} \int_{O} \ \abs{\dfrac{\ratcp{r}{\beta,w}(y)}{\ratcp{r}{\beta,w}(x)}}  \, h(y) \, h(x) \ \mathrm{d}y \, \mathrm{d}x \\
	=& \frac{1}{P[X \in I] \cdot P[Y\in O]} \int_{I} \dfrac{h(x)}{\abs{\ratcp{r}{\beta,w}(x)}}  \, \mathrm{d}x \cdot \int_{O} \abs{\ratcp{r}{\beta,w}(y)} \, h(y) \, \mathrm{d}y.
	\end{align}
\end{subequations}
As $E(\beta,w)$ is non-negative, this directly motivates the minimisation problem\index{Minimisation problem}
\begin{equation}
\label{eqn:minexp}
\argmin{\beta \in \C^m, w \in (\C \setminus \R)^m}{E(\beta,w)}.
\end{equation}
It is possible to compute the value of $E(\beta,w)$ by numerical methods. We have implemented a numerical approach in \tool{MATLAB} and tested different probability distributions of eigenvalues. Two Probabilistic \RFF{}s are depicted in Figure~\ref{fig:exp_conv}, for a gap parameter $G=0.99$. However, they do not yet perform well in practice and require further research on the right choice of parameters and the computation of the occurring integrals. So, we will not focus on them in the analysis of different \RFF{} in the following chapter.

However, if we had computed a solution to \eqref{eqn:minexp}, this would yield both a well-performing \RFF{} and the expected convergence rate.


%% file: figures/experiments/different_lower_bounds_residuals.tex
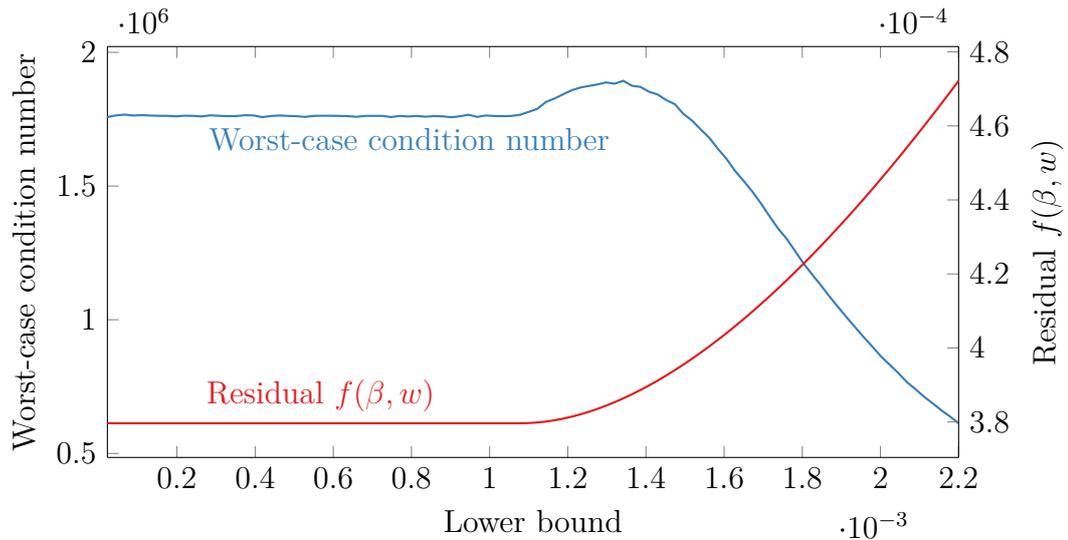
\begin{figure}
	\centering			
	\tikzsetnextfilename{different_lower_bounds_residuals}
	\begin{tikzpicture}
	\begin{axis}[
	scale only axis,
	axis y line*=left,
	xmin=0.000022,
	xmax=0.0022,
	width=0.7\textwidth,
	height=0.24\textheight,
	xlabel={Lower bound},ylabel=Worst-case condition number]
	\addplot[color=GAUSSC, thick] table[x expr=(\thisrowno{0}/100 * 0.0022),y index={2}] {figures/experiments/different_lower_bounds_residuals.dat} node[below, pos=0.05] {Worst-case condition number};
	\end{axis}
	\begin{axis}[
	scale only axis,
	axis y line*=right,
	axis x line=none,
	xmin=0.000022,
	xmax=0.0022,
	ylabel={Residual $f(\beta,w)$},
	width=0.7\textwidth,
	height=0.24\textheight,
	]
	\addplot[color=ZOLOC, thick] table[x expr=(\thisrowno{0}/100 * 0.0022),y index={1}] {figures/experiments/different_lower_bounds_residuals.dat} node[above, pos=0.25] {Residual $f(\beta,w)$};
	\end{axis}
	\end{tikzpicture}
	\caption[Iterations Gamma]{Residuals and worst-case conditions numbers for different lower bounds for the weight function Box-\LSOP/ and $16$ poles.}
	\label{fig:different_lower_bounds_residuals}
\end{figure}

%% file: figures/constraints/nondecreasing.tex
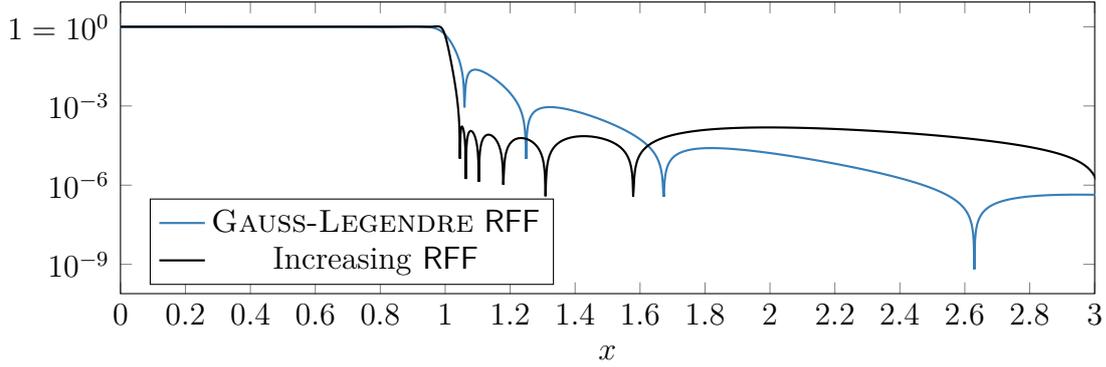
\begin{figure}
	\centering			
	\tikzsetnextfilename{nondecreasing_log}
	\begin{tikzpicture}
	\begin{axis}[
	xmin=0,
	xmax=3,
	ymode=log,
	ytick={0.000000001,0.000001,0.001,1},
	yticklabels={$10^{-9}$,$10^{-6}$,$10^{-3}$,$1 = 10^{0}$},
	legend pos=south west,
	xlabel={$x$}, 
	width=0.9\textwidth,
	height=0.24\textheight,
	]
	\addplot[color=GAUSSC, thick] table[x index={0},y expr=abs(\thisrowno{1})] {figures/RFF/rffs.dat};	
	\addplot[thick] table[x index={0},y expr=abs(\thisrowno{1})] {figures/constraints/nondecreasing.dat};	
	
	\addlegendentry{\textsc{Gauss-Legendre} \RFF{}}
	\addlegendentry{Increasing \RFF{}}
	\end{axis}
	\end{tikzpicture}
	\caption[Absolute value]{Logarithmic absolute function values of non-increasing and increasing $16$-pole \RFF{}s.}
	\label{fig:nondecreasing_log}
\end{figure}

%% file: figures/experiments/exp_conv.tex
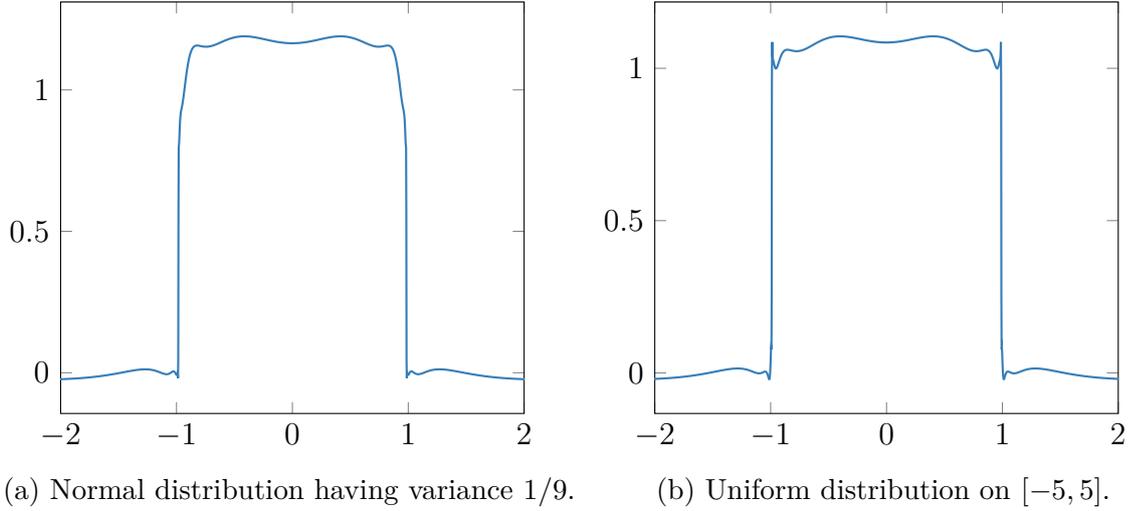
\begin{figure}
    \centering
    \begin{subfigure}[c]{0.48\textwidth}
        \tikzsetnextfilename{expconvgauss}
        \begin{tikzpicture}
        \begin{axis}[
        xmin=-2,
        xmax=2,
        ytick={0,0.5,1},
        width=1\textwidth,
        height=0.31\textheight]
        \addplot[color=GAUSSC, thick] table[x index={0},y index={1}] {figures/experiments/exp_conv_gauss.dat};
        \end{axis}
         \end{tikzpicture}
        \subcaption{Normal distribution having variance $1/9$.}
        \label{fig:exp_conv_gauss} 
    \end{subfigure}
    \begin{subfigure}[c]{0.48\textwidth}
		\tikzsetnextfilename{expconvequi}
		\begin{tikzpicture}
		\begin{axis}[
		xmin=-2,
		xmax=2,
		ytick={0,0.5,1},
		width=1\textwidth,
		height=0.31\textheight]
		\addplot[color=GAUSSC, thick] table[x index={0},y index={1}] {figures/experiments/exp_conv_equi.dat};
		\end{axis}
		\end{tikzpicture}
		\subcaption{Uniform distribution on $[-5,5]$.}
		\label{fig:exp_conv_equi} 
	\end{subfigure}
    \caption[Probabilistic \RFF{}s]{Probabilistic \RFF{}s with $16$ poles.}
    \label{fig:exp_conv}
\end{figure}

%% file: parts/experiments.tex
\preamble{D. Adams}{Last Chance to See}{I [...] am rarely happier than when spending an entire day programming my computer to perform automatically a task that would otherwise take me a good ten seconds to do by hand.}
\chapter{Algorithmic design of weight functions}
\label{chp:tests}

We have previously seen that the \textsc{Gauss-Legendre} \RFF{}\index{RFF!Gauss-Legendre} aims at providing a good behaviour on average, whereas the \textsc{Zolotarev} \RFF{} is for the worst-case convergence scenario. In this chapter, we improve the \textsc{Gauss-Legendre} \RFF{} to ensure fast convergence of \FEAST{} in practice.

Our improvements yield a new \RFF{}, the \textit{Enhanced} $\gamma$-\LSOP/ \RFF{}, which is denoted in Appendix~\ref{chp:resulting_rffs}. This \RFF{} was obtained using the \BFGS{} algorithm, as outlined in Section~\ref{sec:improving_gauss}. In Section~\ref{sec:woc}, we compare the worst-case performance of our new \RFF{} to state-of-art \RFF{}s as introduced in Chapter~\ref{chp:rff}. In Section~\ref{sec:performace}, we examine the performance in practice on a large problem set.

\section{Minimum weight functions}
\label{sec:improving_gauss}

The \LSOP/ minimisation problems aim at minimising the weighted, squared distance to the indicator function. However, this does not necessarily result in an \RFF{} performing well in \FEAST{}, as this process depends on the chosen weight function. In the following, we introduce a criterion for identifying well-performing weight functions. We call this criterion the \textit{worst-case convergence rate of weight functions}\index{Weight function}\index{Worst-case convergence rate!Weight function}.

We will use this worst-case convergence rate of a weight function to improve the $\gamma$-\LSOP/ weight function. This particular weight function was proposed by \textsc{Winkelmann} and \textsc{Di Napoli} to obtain \LSOP/ \RFF{}s\index{SLiSe} to replace the \textsc{Gauss-Legendre} \RFF{}. For $x\in\R$, this weight function is defined as
\begin{equation}
\label{eqn:gamma_slise}
\mathfrak{G}_{\gamma\text{-\LSOP/}}(x) :=
\begin{cases}
1, &\text{if } \abs{x} < 0.95, \\
0.01, &\text{if } 0.95 \le \abs{x} < 1.05, \\
10, &\text{if } 1.05 \le \abs{x} < 1.4, \\
20, &\text{if } 1.4 \le \abs{x} < 5, \\
0, & \text{otherwise.}
\end{cases}
\end{equation}
In particular, we focus on the case of $16$ poles, as this is the standard choice in \FEAST{}. 

We proceed as follows: First, we introduce the worst-case convergence rate of an \RFF{}. Second, we define the worst-case convergence rate of a weight function. Third, we motivate a minimisation problem of the worst-case convergence rate of weight functions. Fourth, we study how to solve the new minimisation problem numerically.
\begin{enumerate}
	\litem{Worst-case convergence rate of RFFs} The worst-case convergence rate has previously been studied for \RFF{}s in \FEAST{}. From the convergence rate of \FEAST{}\index{RFF!Worst-case convergence rate}\index{FEAST}\index{Worst-case convergence rate!RFF} as per \eqref{eqn:conv}, the worst-case convergence rate was derived in \cite{guettel} as \index{Gap!Parameter}\index{Minimisation problem}
	\begin{equation}
		\label{eqn:worst}
		\text{worstcase}(r,G) = \frac{\max_{x \in [G^{-1},\infty]}\abs{r(x)}}{\min_{x \in [0,G]}\abs{r(x)}},
	\end{equation}
	for an \RFF{} $r$ and a gap parameter $G \in (0,1)$, assuming there were no eigenvalues in $[-G^{-1},-G] \cup [G,G^{-1}]$.
	
	As in Definition~\ref{def:zolotarev} of the \textsc{Zolotarev} \RFF, a gap parameter is required for this definition. As before, this is necessary for a meaningful notion of worst-case convergence rate, as \RFF{}s are continuous functions.
	
	\litem{Worst-case convergence rate of weight functions} We now use the worst-case convergence rate of an \RFF{} as per \eqref{eqn:worst}, to assign a worst-case convergence rate to a weight function.
	
	Let $\mathfrak{G}$ be a weight function. Then, we can use $\mathfrak{G}$ in a \LSOP/ minimisation \eqref{eqn:complex_minimiser} with the $16$-pole \textsc{Gauss-Legendre} \RFF{} as the starting point. We choose this \RFF{} as a fixed starting point, as we are interested in improving this \RFF{}. The \LSOP/ minimisation using the weight function $\mathfrak{G}$ yields a new \RFF{} $r$. For this \RFF{} $r$, we have already derived the worst-case convergence rate as per \eqref{eqn:worst}, so we can define
	\begin{equation}
	\text{worstcase}(\mathfrak{G}) := \text{worstcase}(r,G),
	\end{equation}
	where $G \in (0,1)$ is a fixed gap parameter.
	
	The worst-case convergence rate requires a gap. We suggest a moderate choice of $G=0.95$, as this yields good results in practice.
	\litem{Weight functions as a minimisation problem} We we want to improve the previously introduced $\gamma$-\LSOP/ weight function in terms of worst-case convergence rate. To achieve this improvement, we generalise the weight function \eqref{eqn:gamma_slise} to improve to
	\begin{equation}
	\mathfrak{G}^{v,w}(x) :=
	\begin{cases}
	v_1, &\text{if } \abs{x} < w_1, \\
	v_2, &\text{if } w_1 \le \abs{x} < w_1^{-1}, \\
	v_3, &\text{if } w_1^{-1} \le \abs{x} < w_2, \\
	v_4, &\text{if } w_2 \le \abs{x} < w_3, \\
	v_5, & \text{otherwise,}
	\end{cases}
	\end{equation}
	for $x\in\R$, where $v \in (0,\infty)^{5}$, $w_1 \in (0,1)$, $1<w_1^{-1}<w_2<w_3$.\index{Minimisation problem}
	
	The parameters $v_1, \dots, v_5$ denote different positive weights. The values $w_1, w_1^{-1}$ indicate gaps around $\pm 1$. The remaining parameters $w_2,w_3$ offer some more degrees of freedom in the choice of steps of the weight function.
	
	An optimal choice of vectors $v,w$ can be defined as the minimum worst-case convergence rate of the induced weight function $\mathfrak{G}^{v,w}$. This is supposed to yield a better weight function compared to the original $\gamma$-\LSOP/ weight function in \LSOP/. We define the function 
	\begin{equation}
	h(v,w) := \text{worstcase}(\mathfrak{G}^{v,w}),
	\end{equation}
	where $v \in (0,\infty)^{5}$, $w_1 \in (0,1)$, $w_1^{-1}<w_2<w_3$.
	
	This is a \textit{loss function without derivatives} as the convergence rate is always non-negative. This motivates the following minimisation problem of $8$ variables
	\begin{equation}
		\label{eqn:minweight}
		\argmin{v \in (0,\infty)^{5}, \ w_1 \in (0,1), \ w_1^{-1}<w_2<w_3}{h(v,w)}.
	\end{equation}
	\litem{Solving the minimisation problem}
	The minimisation problem \eqref{eqn:minweight} can be solved using standard derivative-free minimisation algorithms. However, the minimisation problem is computationally difficult, because (i) the minimisation problem is non-linear, (ii) the loss function $h$ has no derivatives and (iii) the function evaluations of $h$ are rather expensive. Evaluating $h$ involves a \LSOP/ minimisation and solving the two minimisation problems in \eqref{eqn:worst}. The performance can be increased by using the weights and intervals of the $\gamma$-\LSOP/ weight function as a starting point in the minimisation \eqref{eqn:minweight}. This is a natural choice, as we want to improve the $\gamma$-\LSOP/ weight function in terms of worst-case convergence rate to obtain a replacement candidate for the \textsc{Gauss-Legendre} \RFF{}.
\end{enumerate}

In the end, we obtained a weight function, the \textit{Enhanced} $\gamma$-\LSOP/, which can be found in Example~\ref{exmpl:enhanced_slise}. This new weight function yields a new \RFF{}s via \LSOP/ minimisation. We analyse our new Enhanced $\gamma$-\LSOP/ \RFF{}s in the following. For the standard choice of $16$ poles in \FEAST{}, our new \RFF{} is denoted in Appendix~\ref{chp:resulting_rffs}.

\section{Worst-case performance}
\label{sec:woc}

\input{figures/experiments/convergence_rates}

In the previous section, we obtained a new weight function via minimisation of its worst-case convergence rate. In this section, we analyse the worst-case convergence rates as per \eqref{eqn:worst} of the resulting Enhanced $\gamma$-\LSOP/ \RFF{}s. To achieve this, we compare to its competitors, i.e.\ the \textsc{Gauss-Legendre} and the $\gamma$-\LSOP/ \RFF{}. \textsc{Gauss-Legendre} is the origin of both the other studied \RFF{}s, which were obtained via \LSOP/ minimisation using the \textsc{Gauss-Legendre} \RFF{} as the starting point. The $\gamma$-\LSOP/ \RFF{} was proposed in the original publication as a replacement candidate for the \textsc{Gauss-Legendre} \RFF{} \cite{jan}.

We have computed the worst-convergence rates for different numbers of poles and gap parameters $G$ in Table~\ref{tab:convergence_rates}. In all cases, the \textit{Enhanced} $\gamma$-\LSOP/ provides a better (smaller) worst-case convergence rate than the others.

Recall that we used a gap of $G=0.95$ and $16$ poles when obtaining our new weight function. Surprisingly, our new \RFF{} outperforms its competitors for other gaps and numbers of poles as well.

In the end, we expect our new \RFF{} to perform better in practice. We analyse this in the following section on a large set of benchmark problems.
\section{Performance in practice}
\label{sec:performace}

\input{figures/experiments/iterations_gamma}

In the following, we analyse the performance of our new Enhanced $\gamma$-\LSOP/ \RFF{} in practice. As before, we compare to its competitors, i.e.\ the \textsc{Gauss-Legendre} \RFF{} and the $\gamma$-\LSOP/ \RFF{}.

The performance of \FEAST{}\index{FEAST} in practice has been tested using a set of $2117$ different \HIEP{}s. This is the same set as used by \textsc{Winkelmann} and \textsc{Di Napoli} in \cite{jan}. This set was obtained from the sparse, \textsc{Hermitian} \textit{Si2} matrix, part of the matrix collection of the \textsc{University of Florida} \cite{Si2}. From this matrix, there were selected $2117$ different search intervals, each containing between $5$ and $20$ percent of the eigenvalues of the \tool{Si2} matrix.

\FEAST{}\index{FEAST} was used in the most-recent version \tool{3.0}, compiled using the \tool{Intel Compiler 17.0.0} and and run on an \tool{Intel Core i7-6900K}. We used the \tool{scsrevx} method of \FEAST{} for solving \HIEP{}s of sparse matrices and embedding \RFF{}s. Furthermore, the target accuracy of the eigenvalues was adjusted to $10^{-13}$. In the end, we used the same set-up as in \cite{jan}.\footnote{In fact, in \cite{jan}, the \tool{Intel Compiler} was used in version \tool{16.0.2}. We used the most recent version \tool{17.0.0}.}

As discussed in Section~\ref{sec:feast}, a required argument of \FEAST{} is an upper bound on the number of eigenvalues in the search interval of a given \HIEP{}. Recall that in \FEAST{}, a smaller upper bound yields a better performance in solving the reduced eigenvalue problems, but a worse convergence rate depending on the chosen \RFF{} as per \eqref{eqn:conv}. So, we have tested different upper bounds. For all the studied \HIEP{}s, we determined the actual eigenvalue count in the search interval and scaled it by a factor $C > 1$.  Thus, a smaller $C$ indicates faster \FEAST{} iterations.

The results are depicted in Figure~\ref{fig:iterations_gamma}. For every one of the $2117$ problems, all \RFF{}s assured the convergence of \FEAST{} to the sought eigenvalues. We notice that our new Enhanced $\gamma$-\LSOP/ \RFF{} outperforms the other \RFF{}s, notably for a small $C$. For $C\ge1.1$, the number of required iterations stays almost constant.

\input{figures/experiments/iterations_gamma_1.1_all}

For the case $C=1.1$, we have also compared the actually required counts of \FEAST{} iterations, see Figure~\ref{fig:iterations_gamma_1.1_all}. We note that the \textsc{Gauss-Legendre} \RFF{} terminates after a much larger number of iterations than in the case of the \RFF{}s obtained via weighted, squared approximation. The new Enhanced $\gamma$-\LSOP/ \RFF{} requires only $3$ iterations for most of the \HIEP{}s in contrast to the $\gamma$-\LSOP/ \RFF{}, which needs $4$ iterations in most cases.

We conclude that our new \RFF{}s perform well in \FEAST{} in terms of worst-case and practice.

%% file: figures/experiments/convergence_rates.tex
\begin{table}
	\centering			
	\pgfplotstabletypeset[
	columns={G,q,gauss,bfgsgamma,ebfgs}, 
	columns/G/.style={
		precision=6,
		column name=Gap $G$,
	},
	columns/q/.style={
		int detect,
		column name={Poles},
		preproc/expr={4*##1}
	},		
	columns/{gauss}/.style={
		sci,sci zerofill,precision=2,sci superscript,
		column name={\textsc{Gauss-Legendre}}
	},
	columns/{bfgsgamma}/.style={
		sci,sci zerofill,precision=2,sci superscript,
		column name={$\gamma$-\LSOP/} 
	},
	columns/{box}/.style={
		sci,sci zerofill,precision=2,sci superscript,
		column name={\LBFGSB{}}
	},
	columns/{zolo}/.style={
		sci,sci zerofill,precision=2,sci superscript,
		column name={Elliptic}
	},
	columns/{zolo2}/.style={
		sci,sci zerofill,precision=2,sci superscript,
		column name={Ell. Filter}
	},
	columns/{ebfgs}/.style={
		sci,sci zerofill,precision=2,sci superscript,
		column name={Enhanced $\gamma$-\LSOP/}
	},
	every row no 2/.style={after row=\midrule},
	every row no 5/.style={after row=\midrule},
	every row no 8/.style={after row=\midrule},
	every head row/.style={before row=\toprule,after row=\midrule},
	every last row/.style={after row=\bottomrule},
	col sep=comma
	] %
	{figures/experiments/convergence_rates.dat}
	
	\caption[Convergence rates]{Worst-case convergence rates for different \RFF{}s, gaps and numbers of poles. (smaller is better)}
	\label{tab:convergence_rates}
\end{table}

%% file: figures/experiments/iterations_gamma.tex
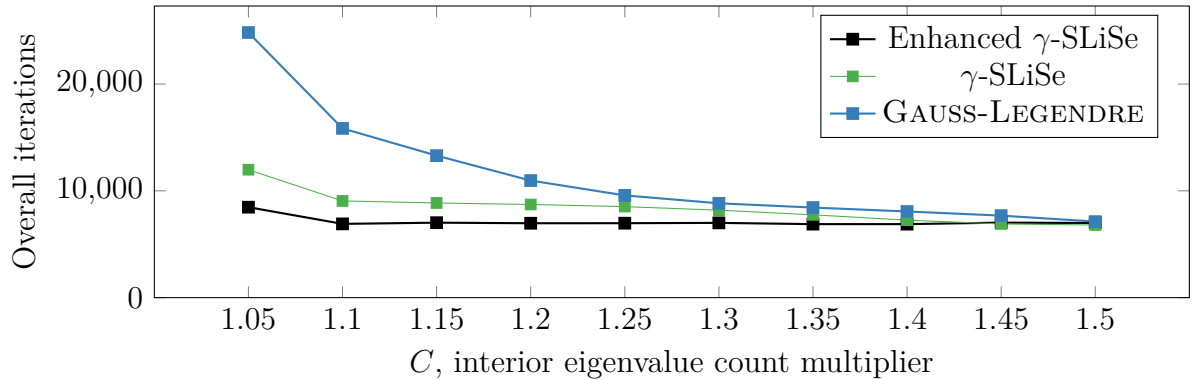
\begin{figure}
	\centering			
	\tikzsetnextfilename{iterations_gamma}
	\begin{tikzpicture}
	\begin{axis}[
	scaled y ticks=false, 
	legend pos=north east,
	ymin=0,
	xmin=1.0,
	xmax=1.55,
	xtick={1.05,1.1,...,1.5,1.5},
	xlabel={$C$, interior eigenvalue count multiplier}, ylabel={Overall iterations},
	width=0.95\textwidth,
	height=0.24\textheight,
	]
	\addplot[mark=square*, thick] table[x index={0},y expr=(\thisrowno{4})] {figures/experiments/iterations_gamma.dat};
	\addplot[color=LSOPC, mark=square*] table[x index={0},y expr=(\thisrowno{2})] {figures/experiments/iterations_gamma.dat};
	\addplot[color=GAUSSC, mark=square*, thick] table[x index={0},y expr=(\thisrowno{1})] {figures/experiments/iterations_gamma.dat};
	\addlegendentry{Enhanced $\gamma$-SLiSe}
	\addlegendentry{$\gamma$-SLiSe}
	\addlegendentry{\textsc{Gauss-Legendre}}
	\end{axis}
	\end{tikzpicture}
	\caption[Iterations Gamma]{Overall number of \FEAST{} iterations required for solving $2117$ \HIEP{}s.}
	\label{fig:iterations_gamma}
\end{figure}

%% file: figures/experiments/iterations_gamma_1.1_all.tex
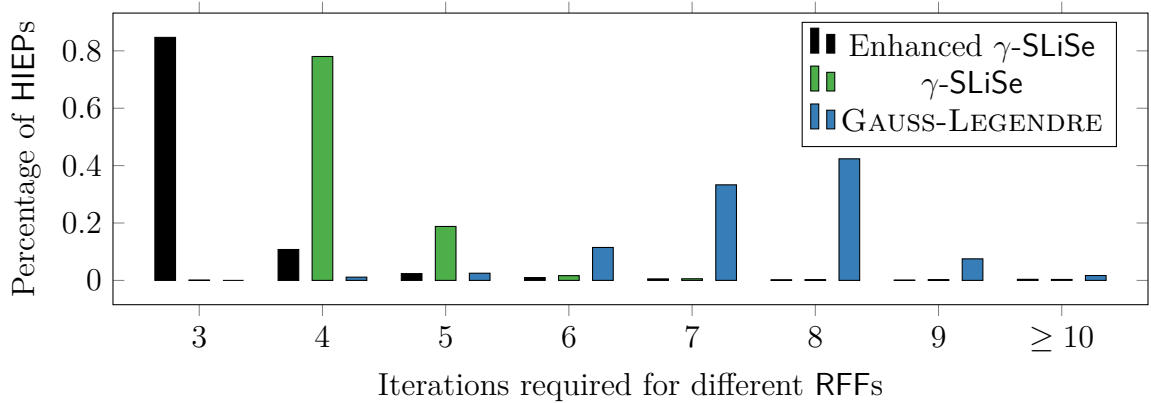
\begin{figure}
	\centering
		\tikzsetnextfilename{iterations_gamma_1.1_all}
	\begin{tikzpicture}
		\begin{axis}[
		xlabel=Iterations required for different \RFF{}s,
		ylabel=Percentage of \HIEP{}s,
		ybar=5pt,
		xticklabels={3,...,9,$\ge 10$},xtick={3,...,10},
		ytick={0,0.2,...,1},
		bar width=0.017\textwidth,
		legend pos=north east,
		width=0.95\textwidth,
		height=0.24\textheight,
		]
		\addplot [fill=black] table[x index={0}, y expr=(\thisrowno{3}/2117)] {figures/experiments/iterations_gamma_1.1_all.dat};
		\addplot [fill=LSOPC] table[x index={0}, y expr=(\thisrowno{4}/2117)] {figures/experiments/iterations_gamma_1.1_all.dat};
		
		\addplot [fill=GAUSSC] table[x index={0},y expr=(\thisrowno{1}/2117)] {figures/experiments/iterations_gamma_1.1_all.dat};
				
		\legend{Enhanced $\gamma$-\LSOP/, $\gamma$-\LSOP/, \textsc{Gauss-Legendre}}
		\end{axis}
	\end{tikzpicture}
	\caption[Counts of iterations]{\FEAST{} iterations for eigenvalue count multiplier $C=1.1$ and 2117 \HIEP{}s.}
	\label{fig:iterations_gamma_1.1_all}
\end{figure}

%% file: parts/conclusion.tex
\chapter{Conclusions}
We have discussed \LSOP/\index{SLiSe} minimisation for designing \RFF{}s, recently introduced in \cite{jan}. In particular, we focussed on current drawbacks of \LSOP/, i.e.
\begin{enumerate}
	\item Incompatibility with standard minimisation algorithms,
	\item Impracticable performance of box-constraints,
	\item Choice of suitable weight functions,
\end{enumerate}
as mentioned in the introductory Section~\ref{sec:intro}.

In this thesis, we have solved these open issues of \LSOP/ and obtained the following results.
\begin{enumerate}
	\litem{Embedding into standard minimisation algorithms} The \LSOP/ minimisation problems have been embedded into standard minimisation algorithms, in particular the \BFGS{} algorithm. This specific algorithm offers the efficiency and extensibility required for \LSOP/. For \LSOP/ minimisation via \BFGS{}, instead of hours, only seconds are needed for obtaining an \RFF{} on an average computer architecture.
	
	\litem{Improving the performance of box-constraints} Box-constraints help reducing the condition numbers in \FEAST{} for \LSOP/ \RFF{}s. We have improved the original implementation of box-constraints in \LSOP/ using the \LBFGSB{} algorithm. Compared to the original implementation of box-constraints in \LSOP/, our approach offers a performance suitable for practice.
	
	We have tested and compared different set-ups of box-constrained \LSOP/. In these set-ups, our \LBFGSB{} approach requires no more than $500$ steps, whereas the original implementation requires up to $900,000$ steps. Furthermore, we have shown that our constrained \LSOP/ \RFF{}s perform significantly better than previous ones, as \LBFGSB{} offers a better accuracy compared to the original implementation.
	
	\litem{Algorithmic choice of weight functions} The choice of weight functions as a parameter of \LSOP/ is crucial for obtaining well-performing \RFF{}s via \LSOP/ minimisation. We defined the worst-case convergence rate of weight functions as a criterion for indicating well-performing weight functions. We derived new \RFF{}s using a weight function with a minimum worst-case convergence rate in \LSOP/.
	
	Our new \RFF{}s offer a better performance in \FEAST{} compared to state-of-art \RFF{}s, such as \textsc{Gauss-Legendre}\index{RFF!Gauss-Legendre} and former \LSOP/ \RFF{}s. Using our new \RFF{}s, the \FEAST{} algorithm requires up to one quarter fewer iterations on average, tested on a large set of benchmark problems.
\end{enumerate}

\textbf{Further work.} For further research, Probabilistic \RFF{}s\index{RFF!Probabilistic} are certainly of major interest. We have presented applications, a mathematical foundation and first, auspicious results of this approach. The numerical computation of Probabilistic \RFF{}s still requires further examination. In the end, we expect both well-performing \RFF{}s and an estimate of the required iterations of \FEAST{} as a result of this approach.

%% file: parts/appendix.tex
\chapter{Ready-to-use rational filter functions}
\label{chp:resulting_rffs}

In the following, we present the studied \RFF{}s, denoted as in \eqref{eqn:rff}.

\begin{center}	
	\footnotesize
	\begin{tabular}{ll}
		\toprule
		\normalsize{Poles \(w\)} & \normalsize{Coefficients \(\beta\)}\\
		\midrule
		\(-0.9980552138505067 + 0.062336105956370486 \iu\)& \(0.02525791710871586 - 0.0015775481910044564 \iu\)\\
		\(-0.9494253842988177 + 0.3139927382100546 \iu \)& \(0.05278354977406013 - 0.017456507483534722 \iu\)\\
		\(-0.7348899387554323 + 0.678186388770961 \iu \)& \(0.05763496444397823 - 0.05318789432545047 \iu\)\\
		\(-0.2841679239019292 + 0.9587745256428074 \iu \)& \(0.025765774438829884 - 0.0869329930919054 \iu\)\\
		\bottomrule
	\end{tabular}
	\captionof{table}{The 16-pole \textsc{Gauss-Legendre} \RFF{}.}
\end{center}

\begin{center}	
	\footnotesize
	\begin{tabular}{ll}
		\toprule
		\normalsize{Poles \(w\)} & \normalsize{Coefficients \(\beta\)}\\
		\midrule
		\(-0.9999975815339606 + 0.0021993013049440135 \iu\)& \(0.0008989201462643977 - 1.977001032029609\text{E-6} \iu\)\\
		\(-0.9998514744807556 + 0.017234528675274002 \iu \)& \(0.005245791227192865 - 9.042216932920706\text{E-5} \iu\)\\
		\(-0.9933358764099828 + 0.11525552757595411 \iu \)& \(0.03462538525214074 - 0.004017540430714314 \iu\)\\
		\(-0.7398348571484926 + 0.6727885136861876 \iu \)& \(0.15051737271560608 - 0.13687697801045523 \iu\)\\
		\bottomrule
	\end{tabular}
	\captionof{table}{The 16-pole \textsc{Zolotarev} \RFF{}.}
\end{center}

\begin{center}	
	\footnotesize
	\begin{tabular}{ll}
		\toprule
		\normalsize{Poles \(w\)} & \normalsize{Coefficients \(\beta\)}\\
		\midrule
		\(-0.9999983713139353 + 0.0022 \iu\)& \(0.0010905705446617412 - 1.4902889756769852\text{E-6} \iu\)\\
		\(-0.9998476756269521 + 0.023174916170735475 \iu \)& \(0.007300520076462563 - 0.00010162408356932002 \iu\)\\
		\(-0.9897979425768154 + 0.15300422557734145 \iu \)& \(0.0435127109551866 - 0.006058193629191226 \iu\)\\
		\(-0.6868662884959791 + 0.7440732728350293 \iu \)& \(0.1355339692180714 - 0.14590122484259907 \iu\)\\
		\bottomrule
	\end{tabular}
	\captionof{table}{The 16-pole \LBFGSB{} Box-\LSOP/ \RFF{}, obtained from and to replace the \textsc{Zolotarev} \RFF{} using a lower bound of $lb=0.0022$.}
\end{center}

\begin{center}	
	\footnotesize
	\begin{tabular}{ll}
		\toprule
		\normalsize{Poles \(w\)} & \normalsize{Coefficients \(\beta\)}\\
		\midrule
		\(-0.9997180876994749 + 0.010064168904151764 \iu\)& \(0.005218903896671892 - 0.0003275342117714203 \iu\)\\
		\(-0.985330269864567 + 0.08344015646402761 \iu \)& \(0.019780578125967584 - 0.005308415315997665 \iu\)\\
		\(-0.8908400599591626 + 0.30261876848986174 \iu \)& \(0.053241710348050676 - 0.03215097589453323 \iu\)\\
		\(-0.43598745582039683 + 0.6982671139969543 \iu \)& \(0.05378661362857605 - 0.12118676200021669 \iu\)\\
		\bottomrule
	\end{tabular}
	\captionof{table}{The 16-pole \BFGS{} \(\gamma\)-\LSOP/ \RFF{}, obtained from and to replace the \textsc{Gauss-Legendre} \RFF{}.}
\end{center}

\begin{center}	
	\footnotesize
	\begin{tabular}{ll}
		\toprule
		\normalsize{Poles \(w\)} & \normalsize{Coefficients \(\beta\)}\\
		\midrule
		\(-0.995102777784057 + 0.01971965034279112 \iu\)& \(0.007451889566376135 - 0.0023538898767857387 \iu\)\\
		\(-0.9656137585011698 + 0.09822459880633161 \iu \)& \(0.019581536492404246 - 0.00823771601370859 \iu\)\\
		\(-0.8531623369434934 + 0.30357032990253513 \iu \)& \(0.04865850681408789 - 0.033809650419106246 \iu\)\\
		\(-0.4113331147792164 + 0.6641012378282691 \iu \)& \(0.04909233881671418 - 0.11480784939181093 \iu\)\\
		\bottomrule
	\end{tabular}
	\captionof{table}{The 16-pole \BFGS{} Enhanced \(\gamma\)-\LSOP/ \RFF{}, obtained from and to replace the \textsc{Gauss-Legendre} \RFF{}.}
\end{center}

\chapter{Considered weight functions}
\label{chp:weight_functions}
\begin{example}[$\gamma$-\LSOP/]
	\label{exmpl:slise}
	\textsc{Winkelmann} and \textsc{Di Napoli} introduced the \textit{$\gamma$-\LSOP/ function} in \cite{jan}, given by
	
	\begin{equation}
	\mathfrak{G}_{\gamma\text{-\LSOP/}}(x) :=
	\begin{cases}
	1, &\text{if } \abs{x} < 0.95, \\
	0.01, &\text{if } 0.95 \le \abs{x} < 1.05, \\
	10, &\text{if } 1.05 \le \abs{x} < 1.4, \\
	20, &\text{if } 1.4 \le \abs{x} < 5, \\
	0, & \text{otherwise,}
	\end{cases}
	\end{equation}
	for all $x \in \R$. The minimiser of the resulting minimisation problem is to replace the \textsc{Gauss-Legendre} \RFF{} as of Section~\ref{sec:gauss}.
\end{example}

\begin{example}[Box-\LSOP/]
	\label{exmpl:box_slise}
	\textsc{Winkelmann} and \textsc{Di Napoli} introduced the \textit{Box-\LSOP/ function} in \cite{jan}, given by
	
	\begin{equation}
	\mathfrak{G}_{\text{Box-\LSOP/}}(x) :=
	\begin{cases}
	1, &\text{if } \abs{x} < 0.95, \\
	4, &\text{if } 0.95 \le \abs{x} < 0.995, \\
	2, &\text{if } 0.995 \le \abs{x} < 1.005, \\
	4, &\text{if } 1.005 \le \abs{x} < 1.05, \\
	0.6, &\text{if } 1.05 \le \abs{x} < 1.1, \\
	1, &\text{if } 1.1 \le \abs{x} < 1.3, \\
	0.3, &\text{if } 1.3 \le \abs{x} < 1.8, \\
	0.1, &\text{if } 1.8 \le \abs{x} < 3, \\
	0, & \text{otherwise,}
	\end{cases}
	\end{equation}
	for all $x \in \R$. The minimiser of the resulting minimisation problem is to replace the \textsc{Zolotarev} \RFF{}s as of Section~\ref{sec:zolotarev} when subject to box-constraints.
\end{example}

\begin{example}[Enhanced $\gamma$-\LSOP/]
	\label{exmpl:enhanced_slise}
	For the analysis, we have come up with a new weight function as per
	
	\begin{equation}
	\mathfrak{G}_{\text{Enhanced }\gamma\text{-\LSOP/}}(x) :=
	\begin{cases}
	0.7, &\text{if } \abs{x} < 0.96, \\
	0.00092, &\text{if } 0.96 \le \abs{x} < 1.0417, \\
	887, &\text{if } 1.0417 \le \abs{x} < 1.4, \\
	20, &\text{if } 1.4 \le \abs{x} < 10, \\
	0, & \text{otherwise,}
	\end{cases}
	\end{equation}
	for all $x \in \R$. The minimiser of the resulting minimisation problem is to replace the \textsc{Gauss-Legendre} \RFF{} as of Section~\ref{sec:gauss}.
\end{example}

\input{parts/BFGS}

\chapter{Further details}
\begin{lemma}	
	\label{lem:indicator_integral}
	Let $\gamma: [0,1] \rightarrow \C$ be a closed, positively-oriented, non-self-intersecting path in $\mathbb{C}$ enclosing $(a,b)$. Then, it holds that
	\begin{equation}
	i_{(a,b)}(x) =  \frac{1}{2 \pi \iu} \, \int_{\gamma} \! \frac{1}{\zeta - x} \, \mathrm{d}\zeta, \quad x \in \R \setminus \Ima(\gamma).
	\end{equation}
	
	\begin{proof}
		Fix $z \in \C \setminus \Ima(\gamma)$. From the \textsc{Laurent} expansion, it follows that $\zeta \mapsto \frac{1}{\zeta-z}, \quad \zeta \in \C \setminus \set{z}$ has one residue of $1$ at $z$. Then, the residue theorem implies the result for complex $z$ and thus for real ones as well.
	\end{proof}
\end{lemma}

%

\begin{proposition}
    \label{theorem:linesearch}
    Assume the preliminaries of Definition~\ref{def:wolfe}. Then, there always exists a non-empty interval $I \subseteq (0,\infty)$ such that any $\alpha \in I$ satisfies the \textsc{Wolfe} conditions.
    \begin{proof}
        Define $\Phi (\alpha) := f(x - \alpha p) - (f(x) - c_1 \alpha \nabla \trans{f} p)$, for $\alpha \in \R$. Thus, we have that $\Phi'(\alpha) = (\trans{\nabla f(x + \alpha p) - c_1  \nabla f (x))} p$, for $\alpha \in \R$. As $p$ is a descent direction, we have $\Phi'(0) < 0$. Now, there exists an $\alpha_1 > 0$ such that $\Phi'(\alpha)<0$ for $\alpha \in [0,\alpha_1)$ as $\Phi$ is continuously differentiable. Hence, $\Phi(\alpha_1) < 0$ as well. Then, as $f(x) \ge 0$ for $x \in \R^n$, it is
        \begin{equation}
        \Phi(\alpha) \ge- (f(x) - c_1 \alpha \nabla \trans{f} p) \ge -f(x).
        \end{equation}
        For $\alpha \rightarrow \infty$, the middle term diverges to $\infty$ and so does $\Phi$. Thus, we can find an $\alpha_2 > \alpha_1$ such that $\Phi(\alpha_2) = 0$ by the intermediate value theorem. Choose the smallest such $\alpha_2$. We have that $\Phi(\alpha) < 0$ and hence $\alpha$ satisfies the first \textsc{Wolfe} condition~\eqref{eqn:wolfe1} for $0< \alpha < \alpha_2$.
        
        Now, by Rolle's theorem, there exists an $a_3 \in (0, \alpha_2)$ such that $\Phi'(\alpha_3) = 0$, i.e.
        \begin{equation}
        \trans{\nabla f(x + \alpha_3 p)}p = c_1 \trans{\nabla f(x)} p > c_2 \trans{\nabla f(x)} p,
        \end{equation}
        in other words, the second \textsc{Wolfe} condition~\eqref{eqn:wolfe2}, as $c_1 < c_2$. The continuity of both $f$ and $\nabla f$ ensures that there is an entire interval $I \subseteq (a_3,a_2)$, which proves the claim.
    \end{proof}
\end{proposition}

\begin{theorem}
    \label{thm:bfgs_pd}
    Let $k \in \N_0$, $f$ be a loss function of $n \in \N$ variables and $H_k \in \R^{n \times n}$ symmetric positive-definite. Let $x_k$ be the current iterate, $p_k := -H_k \nabla f (x_k)$ the descent direction in a line search $\alpha \mapsto f(x_k + \alpha p_k)$ and $\alpha_k \in (0,1)$ chosen to satisfy the \textsc{Wolfe} conditions~\eqref{eqn:wolfe} for some parameters $0 < c_1 < c_2 < 1$. Then the next iterate is given by $x_{k+1} := x_k + \alpha p_k$. Define $s_k := x_{k+1} - x_k$ and $y_k := \nabla f(x_{k+1}) - \nabla f(x_{k})$, then the \BFGS{} update $H_{k+1}$ of $H_k$ as per \eqref{eqn:bfgs} is symmetric positive-definite.
    \begin{proof}
        As the symmetric matrices form a group, from \eqref{eqn:bfgs} the symmetry of $H_{k+1}$ follows immediately. By the second \textsc{Wolfe} condition~\eqref{eqn:wolfe2} and as $p_k$ is a descent direction we find
        \begin{equation}
        \nabla f(x_k)^T p_k < c_2 \underbrace{\nabla f(x_k)^T p_k}_{<0} \overset{\mathrm{\eqref{eqn:wolfe2}}}{\le} \nabla f(x_{k+1})^T p_k.
        \end{equation}
        By the definition of $s_k$, it follows that
        \begin{equation}
        \nabla f(x_k)^T s_k < \nabla f(x_{k+1})^T s_k.
        \end{equation}
        Thus, it is
        \begin{equation}
        \label{eqn:spd_helper2}
        y_k^T s_k > 0.
        \end{equation}
        By definition of $H_{k+1}$, and as $H_k$ is positive-definite, we have
        \begin{equation}
	        \underbrace{v^T H_k v}_{\ge 0} + (s_k^{T} z)^2 \underbrace{(y_k^T s_k)^{-1}}_{> 0, \text{ by } \eqref{eqn:spd_helper2}} \ge 0,
        \end{equation}
        where
        \begin{equation}
	        \label{eqn:spd_helper}
	        v := z - y_k \frac{s_k^{T} z}{y_k^T s_k}.
        \end{equation}        
        This implies that $H_{k+1}$ is positive definite and the claim holds.
    \end{proof}
\end{theorem}

\begin{theorem}
	\label{thm:complexcostfunction}
	Let $\mathfrak{G}$ a piece-wise constant weight function. Furthermore, let $\beta \in \C^m, w \in (\C \setminus \R)^m$. For the function $f$ in \eqref{eqn:complexlossfunction}, it then holds that
	\begin{equation}
		\label{eqn:loss}
		f(\beta, w) = 2 \Re[{\trans{\beta} ((X-Z) \beta + (W-Y) \conj{\beta} - 2\theta )}] + \frac{1}{2} \int_{-\infty}^{\infty} \! \mathfrak{G} (x) \, |h(x)|^2 \mathrm{d}x
	\end{equation}
	and for the gradient,
	\begin{subequations}
		\begin{equation}
			\nabla_\beta f(\beta, w) = 4 \, [\ctrans{\beta} (\nabla X - \nabla Z) + \trans{\beta} (\nabla \conj{W} - \nabla \conj{Y}) - \ctrans{\nabla \theta}] \, I_\beta,
		\end{equation}
		\begin{equation}
			\nabla_w f(\beta, w) =  4 \, [\ctrans{\beta} (X - Z) + \trans{\beta} (\conj{W} - \conj{Y}) - \ctrans{\theta}],
		\end{equation}
	\end{subequations}
	where
	\begin{equation*}
		\begin{aligned}[l]
			W_{k,l} &:= \int_{-\infty}^{\infty} \! \frac{\mathfrak{G} (x)}{(x-\conj{w_k})(x-\conj{w_l})} \mathrm{d}x \\
			X_{k,l} &:= \int_{-\infty}^{\infty} \! \frac{\mathfrak{G} (x)}{(x-\conj{w_k})(x-{w_l})} \mathrm{d}x\\
			Y_{k,l} &:= \int_{-\infty}^{\infty} \! \frac{\mathfrak{G} (x)}{(x-\conj{w_k})(x+\conj{w_l})} \mathrm{d}x\\
			Z_{k,l} &:= \int_{-\infty}^{\infty} \! \frac{\mathfrak{G} (x)}{(x-\conj{w_k})(x+{w_l})} \mathrm{d}x\\
			\theta_{k}&:= \int_{-\infty}^{\infty} \! \frac{\mathfrak{G} (x) i_\gamma(x)}{(x-\conj{w_k})} \mathrm{d}x\\
		\end{aligned}
		\qquad
		\begin{aligned}[r]
			\nabla \conj{W_{k,l}} &:= \int_{-\infty}^{\infty} \! \frac{\mathfrak{G} (x)}{(x-{w_k})(x-{w_l})^2} \mathrm{d}x \\
			\nabla X_{k,l} &:= \int_{-\infty}^{\infty} \! \frac{\mathfrak{G} (x)}{(x-\conj{w_k})(x-{w_l})^2} \mathrm{d}x\\
			\nabla \conj{Y_{k,l}} &:= \int_{-\infty}^{\infty} \! \frac{\mathfrak{G} (x)}{(x-{w_k})(x+{w_l})^2} \mathrm{d}x\\
			\nabla Z_{k,l} &:= \int_{-\infty}^{\infty} \! \frac{\mathfrak{G} (x)}{(x-\conj{w_k})(x+{w_l})^2} \mathrm{d}x\\
			\nabla \conj{\theta_{k}} &:= \int_{-\infty}^{\infty} \! \frac{\mathfrak{G} (x) i_\gamma(x)}{(x-{w_k})^2} \mathrm{d}x\\
		\end{aligned}
	\end{equation*}
	are standard definite integrals, for $k,l = 1, \dots, m$.
	\begin{proof}
		Confer \cite{jan}.
	\end{proof}
\end{theorem}

%% file: parts/BFGS.tex
\chapter{Details of BFGS minimisation}
\label{chp:bfgs}

\input{figures/BFGS/scheme}

In Section~\ref{chp:nonlin}, we have seen how to describe \enquote{good} approximations of the indicator function by \RFF{}s as per \eqref{eqn:rff} in terms of a loss function as per Definition~\ref{def:loss}. We expect that minimisers of this loss function yield fast convergence in eigenvalue computation. We present how to find minimisers of such functions efficiently, employing the \BFGS{} minimisation algorithm. We first introduce a general minimisation scheme and then present the \BFGS{} algorithm.

To minimise a loss function $f:\R^n \rightarrow [0,\infty)$ of $n \in \N$ variables in terms of Definition~\ref{def:loss}, for $k=0,1,2,\dots$, one approximates a minimiser of $f$ by \textit{iterates} $x_k \in \R^n$, such that
\begin{equation}
	f(x_0) > f(x_1) > \dots > f(x_k) \ge 0.
\end{equation}

The only question is how to obtain these iterates. We proceed as shown in Figure~\ref{fig:optimisation_scheme}. Given a starting iterate $x_k$, we want to find a next iterate $x_{k+1}$ such that $f(x_k) > f(x_{k+1})$. This yields a better approximation to a local minimiser. However, this equivalent to finding a \textit{descent direction} $p_k \in \R^n$. A descent direction is a vector that \enquote{points} from $x_k$ into a direction where $f$ attains a smaller value. It can be shown by fundamental methods of calculus, i.e. \textit{first-order \textsc{Taylor} expansion}, that this is equivalent to
\begin{equation}
	\label{eqn:descent_direction2}
	{\nabla f(x_k)^T p_k < 0}.
\end{equation}
This vector $p_k$ only shows the direction to a smaller function value of $f$ at $x_k$. In order to gather a smaller function value and thus a next iterate $x_{k+1}$, we have to move along the vector $p_k$. This process is called \textit{line search} as we search along the line
\begin{equation}
	\label{eqn:linesearch}
	\alpha \mapsto f(x_k + \alpha p_k), \quad \text{for} \ \alpha \in (0,\infty)
\end{equation}
to gather an $\alpha_k \in (0,\infty)$ such that $f(x_k) > f(x_{k} + \alpha_k p_k)$. An $\alpha_k$ is called \textit{step length}. In the end, we set $x_{k+1}:=x_{k} + \alpha_k p_k$ and repeat the process for $x_{k+1}$ instead of $x_k$.

\begin{example}
	\label{exmpl:descent}
	A simple example of a descent direction at $x_k$ is given by the negative gradient $-\nabla f(x_k)$. Obviously, this satisfies \eqref{eqn:descent_direction2}, if and only if $\nabla f(x_k) \neq 0$. However, this constraint is no problem. If we have $\nabla f(x_k) = 0$, we have already found (local) \textit{minimiser} if it is no \textit{saddle points} of $f$. In practical implementations, this assured by the right choice of the starting iterate $x_0$.
\end{example}

\begin{example}	
	To illustrate a descent direction, consider the following \tool{2D} case: If we want to find the minimiser of a convex function, e.g. $f(x) := x^2$, we can simply \enquote{follow} the gradient $-\nabla f(x) = -2x$: For any $x \in \R \setminus \set{0}$, this indicates the direction to the minimiser of $f$.
\end{example}

\input{figures/BFGS/contour}

\begin{example}
	\label{exmpl:loss}
	As an example for minimisation of a loss function, we study the function $f(x,y) := x^2 + y^2$ for $x,y \in \R$ as in Example~\ref{exmpl:loss_function} and Figure~\ref{fig:quadratic1}.
	
	All the following steps are depicted in Figure~\ref{fig:contour}. The gradient of $f$ is $\nabla f(x,y) = (2x,2y)$. Suppose, our current iterate was $x_k = (1,1)$. Then, to gather a next iterate, we will compute a descent direction at $x_k$. For this, we are going to use the negative gradient of $f$ at $x_k$ as in Example~\ref{exmpl:descent}. We have for the descent direction $p_k := -\nabla f(x_k) = (-2,-2)$. We cannot just use $x_k + p_k$ as the next iterate as we have $f(x_k + p_k) = f(-1,-1) = f(1,1)$, but we are searching for a decrease in function value. For instance, we choose $\alpha_k = 0.25$. Then, we set $x_{k+1} := x_k + \alpha_k p_k = (0.5,0.5)$.
\end{example}

\section*{Line search}
In Example~\ref{exmpl:descent}, we have shown how to obtain a descent direction. Following Figure~\ref{fig:optimisation_scheme}, we have to perform a line search to gather a new iterate. In this section, we study the line search as per \eqref{eqn:linesearch} more deeply. We note that for an arbitrary choice of step lengths $\alpha_k$ in \eqref{eqn:linesearch}, this does not necessarily yield a minimiser.

\begin{example}
	\label{exmp:non-converging_line_search}
	As an example of a non-converging line search, consider the loss function $f(x) := x^2/2$, for $x \in \R$. The gradient of $f$ is given by $\nabla f(x) = x$. To obtain a minimiser, we choose $x_0 :=-3$ as the starting point.

	If we choose $x_{k+1}$ such that $x_{k+1} - x_k = 2^{-k}$ in every iteration, we have
	\begin{equation}
		x_{k+1} = x_0 + \sum_{n=0}^{k} 2^{-n}
	\end{equation}
	We remember from calculus, that
	\begin{equation}
	\sum_{k=0}^{\infty} 2^{-k} = 2.
	\end{equation}
	In other words,
	\begin{equation}
		\lim_{k \rightarrow \infty} x_{k} = -1 \neq 0.
	\end{equation}
	Hence, we never reach the minimiser of $f$.
\end{example}

If the decrease in $f$ is not sufficiently large, we may never reach a minimiser (nor even get close to it).

\input{figures/BFGS/sufficient_decrease}

To solve this issue, in 1969, Philip Wolfe combined two conditions, that help avoiding no or slow convergence during line search \cite{wolfe1,wolfe2}: We have seen in the previous Example~\ref{exmp:non-converging_line_search} that decrease in function value cannot be arbitrarily. The first \textsc{Wolfe} condition states that there must be so-called \textit{sufficient decrease}. To phrase this, note that for a descent direction $p_k$ the term $\nabla f(x_k)^T p_k$ indicates the \enquote{steepness} of the descent at $x_k$, i.e. the \textit{directional derivative}. We now want for the function value at the next iterate $x_{k+1}$ to lie below the slope $c_1 \alpha_k \nabla f(x_k)^T$ where $c_1 \in (0,1)$ is a control parameter. This imposes on the step length $\alpha_k$ and thus the next iterate $x_{k+1} = x_k + \alpha_k p_k$ the condition that
\begin{equation}
	\label{eqn:first_wolfe}
	f(x_{k+1}) \le f(x_k) + c_1 \alpha_k \nabla f(x_k)^T p_k.
\end{equation}
The condition is also known as the \textit{\textsc{Armijo} rule}, introduced in 1966 \cite{armijo}. Per se, this condition is not sufficient as any step length small enough satisfies the condition.

\begin{example}
	An example is depicted in Figure~\ref{fig:sufficient_decrease}. There, the function and the slope, that shall bound the function, are shown. We find that for $\alpha \in [0,1] \cup [2,2.5]$ the first \textsc{Wolfe} condition as per \eqref{eqn:first_wolfe} is fulfilled.
\end{example}

We know that any minimiser of $f$ satisfies $\nabla f(x) = 0$. Thus, we would like to ensure sufficient decrease in the gradient as well. This second \textsc{Wolfe} condition is known as \textit{curvature condition} and helps during minimisation to reach a minimiser potentially faster. For instance, we would like for the descent direction $p_k$ at the next iterate $x_{k+1}$, that it is at least $c_2  \in (0,1)$ percent as steep as at the current iterate $x_k$. This can be phrased as
\begin{equation}
	c_2 \nabla f(x_k)^T p_k \le \nabla f(x_{k+1})^T p_k \le 0
\end{equation}

\begin{definition}[\textsc{Wolfe} conditions]
	\label{def:wolfe}
	Let $n \in \N$ and $f$ a loss function of $n$ variables. Furthermore let $x \in \R^n$ and $p$ be a descent direction of $f$ at $x$, i.e. $\nabla \trans{f(x)} p < 0$. Then, a step length $\alpha \in (0,\infty)$ of a line search as per \eqref{eqn:linesearch} is said to satisfy the \textit{Wolfe conditions}, if
	\begin{subequations}
		\label{eqn:wolfe}
		\begin{align}	
		\label{eqn:wolfe1}
		f(x + \alpha p) \, \le& \, f(x) + c_1 \alpha \nabla \trans{f(x)} p, \\
		\label{eqn:wolfe2}
		c_2 \nabla f(x)^{T} p \, \le& \, \nabla f(x + \alpha p)^{T} p,
		\end{align}
	\end{subequations}
	for fixed parameters $0 < c_1 < c_2 < 1$.
\end{definition}

According to Nocedal and Wright \cite{numericaloptmization}, a good choice for the mentioned parameters is $c_1 = 10^{-4}$ and $c_2=0.9$.

Now, the Definition~\ref{def:wolfe} of the \textsc{Wolfe} condition raises two questions: Does there always exist a step length satisfying the \textsc{Wolfe} conditions and if so, is this step length computable with a finite number of steps? In fact, this is true. For details, see Proposition~\ref{theorem:linesearch}. Basically, the Proposition states that there always exists an interval $I \subseteq (0,\infty)$ of step lengths satisfying the \textsc{Wolfe} conditions. Such an interval can be found in finitely many steps.

There exists a variety of algorithms computing such a step size, confer \cite[60--62]{numericaloptmization}, \cite{moerthuente}. Especially, the latter presents a popular and efficient implementation. There are also other approaches to line search that perform very well. One of such is the method of Hager and Zhang \cite{hagerzhang} that does not always assure first \textsc{Wolfe} condition~\eqref{eqn:wolfe1} and achieves good results in practice.

\section*{BFGS algorithm}
Now that we know how to find \enquote{good} step lengths in a line search as in \eqref{eqn:linesearch}, we would like to find out how to gather \enquote{good} descent directions. Recall, that, given a current iterate $x_k$, our goal is to find to next iterate $x_{k+1}$ such that $f(x_k) > f(x_{k+1})$.

In minimisation, in each iteration, one usually employs an approximation of $f$, called \textit{model}. We do not need a model of the entire $f$, but only for values \enquote{close} to the current iterate $x_k$. In the end, our hope is that this model will be easier to minimise than $f$. Even more, we believe that a minimiser of the model will yield the direction to a minimiser of $f$. We have seen before that this is called a descent direction.

Say, the minimiser of the model was a variable called $x_\text{min}$. Then, the descent direction to analyse is $p_k := x_\text{min} - x_k$. Recall Figure~\ref{fig:optimisation_scheme}. We can now do a line search for the next iterate along $p_k$. We have to do that as $x_\text{min}$ does not need to be a minimiser of $f$. We have to verify that. Even more, we would like to ensure the \textsc{Wolfe} conditions as introduced in the previous section in Definition~\ref{def:wolfe}.

For the model, it is a straightforward idea to use a \textit{second-order \textsc{Taylor} approximation} of $f$ at $x_k$. This is a quadratic approximation of $f$ in the neighbourhood of $x_k$ and is given by
\begin{subequations}
	\label{eqn:newton}
	\begin{equation}
	f(x) \approx \tilde{f}(x) := f(x_k) + \nabla f(x_k)^T x + \frac{1}{2} x^T \nabla^2 f(x_k) x, \quad x \in \R^n,
	\end{equation}
	at each iteration $k$.
	
	We assume the \textsc{Hessian} $\nabla^2 f(x_k)$ to be symmetric postive-definite. We make this assumption as any positive-definite matrix is non-singular. It follows that the unique minimiser of $\tilde{f}$ is given by
	\begin{equation}
	\label{eqn:newton_les}
	\tilde{p}_{k} = - \nabla^2 f(x_k)^{-1} \nabla f(x_k).
	\end{equation}
\end{subequations}

This approach is called \textit{Newton's method} \cite{numericaloptmization} and usually gives a good approximation of $f$. The computational costs, i.e. solving the linear system~\eqref{eqn:newton_les}, are quite high. In the case of our eigenvalue minimisation, the \textsc{Hessian}s usually have a large condition number. This means that numerically solving \eqref{eqn:newton_les} can only be accomplished with little accuracy.

	\input{figures/BFGS/BFGS_scheme}

There exist multiple efficient approaches to replacing the \textsc{Hessian}s $\nabla^2(x_k)$ by approximations $B_k$. This means that we use the model
\begin{subequations}
	\begin{equation}
	\label{eqn:bfgs_model}
	f(x) \approx m_k(x) := f(x_k) + \nabla f(x_k)^T x + \frac{1}{2} x^T B_k x, \quad x \in \R^n.
	\end{equation}	
	This is called a \textit{quasi}-\textsc{Newton} minimisation. The \BFGS{} algorithm is a powerful example of such. The key feature of this method is the preservation of the positive-definiteness of a starting \textsc{Hessian} approximation $B_0$. In general, the identity matrix $I_n$ is a sufficiently good $B_0$. In each iteration, we can again find the unique minimiser of \eqref{eqn:bfgs_model} as
	\begin{equation}
	\label{eqn:descent_direction}
	p_{k} = - B_{k}^{-1} \nabla f(x_k),
	\end{equation}
	because the positive-definite matrix $B_k$ is non-singular. A scheme of the resulting \BFGS{} algorithm and its required input and produced output is depicted in Figure~\ref{fig:BFGS_scheme}.
	
	Note that a positive-definite $B_{k}$ yields
	\begin{equation}
	\nabla f(x_k)^{T} p_k = - \nabla f(x_k)^{T} B_{k}^{-1} \nabla f(x_k) < 0.
	\end{equation}
	Thus, $p_k$ is a descent direction as per \eqref{eqn:descent_direction2}, which is necessary to be able to obtain a new iterate $x_{k+1}$ via line search.
	
	It seems like solving a linear system was needed to obtain this descent direction $p_k$. In fact, it is possible to compute the inverse of $B_k$ efficiently and analytically. We will denote this inverse as $H_k = B_k^{-1}$. In the end, \eqref{eqn:descent_direction} becomes
	\begin{equation}
		\label{eqn:descent_direction3}
		p_{k} = - H_k \nabla f(x_k).
	\end{equation}
\end{subequations}
We now know how to compute the next descent direction from a given inverse \textsc{Hessian} approximation $H_k = B_k^{-1}$ and the current iterate $x_k$, see also Figure~\ref{fig:BFGS_scheme}. It remains to show how to obtain the next inverse \textsc{Hessian} approximation.

\input{figures/BFGS/BFGS_scheme_final}

For the very first iteration of the \BFGS{} algorithm, i.e. $k=0$, we just assume a symmetric positive-definite $H_0$ as an input, e.g. the identity matrix. From this, we can gather the very first descent direction $p_k$ easily as of the previous equation \eqref{eqn:descent_direction3}. This then yields the next iterate $x_{k+1}$ via line search.

We now have to acquire a next inverse \textsc{Hessian} approximation $H_{k+1}$. For this, we use the previously gathered results to update $H_k$ to a new $H_{k+1}$. This is the \textit{\BFGS{} update} step and given by
\begin{equation}
	\label{eqn:bfgs}
	H_{k+1} = (I_n - \frac{s_k y_k^{T}}{y_k^T s_k}) H_k (I_n - \frac{y_k s_k^{T}}{y_k^T s_k}) +  \frac{s_k s_k^{T}}{y_k^T s_k}
\end{equation}
where $s_k := x_{k+1} - x_k$ and $y_k := \nabla f(x_{k+1}) - \nabla f(x_{k})$.

It can be shown that the \BFGS{} update maintains the symmetric positive-definiteness of a starting inverse \textsc{Hessian} approximation $H_0$, see Theorem~\ref{thm:bfgs_pd}.

We can update Figure~\ref{fig:BFGS_scheme} to show the entire minimisation procedure for the \BFGS{} algorithm using and producing the previously mentioned inputs and outputs in Figure~\ref{fig:BFGS_scheme}, see Figure~\ref{fig:BFGS_scheme_final}. Eventually, we gather the following algorithm
\begin{algorithm}[H]
	\caption{(Unconstrained \BFGS{} Algorithm).}
	\label{algo:bfgs}
	\textbf{Input:} Loss function $f$ of $n \in \N$ variables and its gradient $\nabla f$, $x \in \R^n$, symmetric positive-definite $H \in \R^{n \times n}$, $\eps >0$.\\
	\textbf{Output:} $x \in \R^n$ such that $f(x) \le \eps$.
	\begin{algorithmic}[1]
		\Function{BFGS}{$f, \nabla f, x, H, \eps$}			
			\While{$f(x) > \eps$}
				\State $p \gets -H \nabla f(x)$ \Comment{Compute descent direction}
				\State \label{algo:bfgs_linesearch} $\alpha \gets$ Step length of a line search~\eqref{eqn:linesearch} along $p$ satisfying~\eqref{eqn:wolfe}
				\State $x_{pre} \gets x$
				\State $s \gets \alpha \, p$
				\State $x \gets x_{pre} + s$
				\State $y \gets \nabla f(x) - \nabla f(x_{pre})$
				\State $H \gets$ \BFGS{} update of $H$ according to~\eqref{eqn:bfgs}
			\EndWhile
			\State \textbf{return} $x$
		\EndFunction
	\end{algorithmic}
\end{algorithm}

\begin{observation}
	Let $f$ be a loss function of $n \in \N$ variables and assume that the evaluation of $\nabla f$ and $\nabla^2 f$ can be achieved in constant time and the line search in $\mathcal{O}(g(n))$ for some polynomial $g$. Then, one iteration of the \BFGS{} Algorithm~\ref{algo:bfgs} has a time complexity of $\mathcal{O}(n^2 + g(n))$ in contrast to $\mathcal{O}(n^3 + g(n))$ in \textsc{Newton}'s method~\eqref{eqn:newton}. However, possibly, \textsc{Newton}'s method requires less iterations as no \textsc{Hessian} approximation is used.
\end{observation}

\section*{Handling box-constraints}
\label{sec:constrained_optimisation}

\input{figures/BFGS/box}

In Section~\ref{sec:box}, we have seen how to improve \LSOP/ \RFF{}s by box-constraints. In a general real minimisation environment, box-constraints are known as demanding arbitrary upper and lower bounds on the minimiser. For the $i$-th component of an $x\in\R^n$, we would like to have
\begin{equation}
l_i \le x_i \le u_i,
\end{equation}
where $l_i, u_i \in [-\infty, \infty]$. Here, $l_i$ denotes a lower and $u_i$ an upper bound. This means, that for a loss function $f$ of $n \in \N$ variables, we would like to solve the minimisation problem
\begin{subequations}
	\label{eqn:constrained_optimisation}
	\begin{equation}
	\argmin{x \in \R^n}{f(x)}
	\end{equation}
	\begin{equation}
	\text{subject to:} \ l_i \le x_i \le u_i, \quad \text{for} \ i=1,\dots,n.
	\end{equation}
\end{subequations}
The lower and upper bounds are called \textit{box-constraints}. It is called \enquote{box} as indeed the values in the bounds lie in a box of the common sense.
\begin{example}	
	An example of a box is depicted in Figure~\ref{fig:box}. Here, the box of values lying in the bounds is denoted as $C$.
\end{example}
\begin{example}
	Suppose, we want to impose constraints on an minimisation problem in $\R^2$. For instance, we would like to have for the first dimension that $l_1=2\le x_1 \le8=u_1$ and for the second dimension that $l_2 = -5 \le x_2 \le \infty = u_2$. Thus, all these values lie in a set given by
	\begin{equation}
	C=[2,8] \times [-5,\infty).
	\end{equation}
\end{example}

For the solution of box-constrained problems, there is a very well-known algorithm, called \LBFGSB{} \cite{lbfgs1}. In terms of the \BFGS{} scheme in Figure~\ref{fig:BFGS_scheme_final}, mainly, the used model is altered to incorporate box-constraints yielding the method shown in Figure~\ref{fig:BFGS_scheme_contraint}. We will shortly describe the main idea behind this algorithm:

\input{figures/BFGS/scheme_constraint}

At the current iterate $x_k$, similar to \eqref{eqn:bfgs_model}, the algorithm forms a quadratic model in the neighbourhood of $x_k$ using a \textsc{Hessian} approximation $B_k$ yielding some function
\begin{equation}
\label{eqn:quadratic_model}
q_k(x), \quad \text{for} \ \in \R^n.
\end{equation}
Consequently, the algorithm identifies the components of $x_k$ that are going to violate the bounds in the next iterate $x_{k+1}$ as follows. A \textit{projected line search} along the gradient $\nabla f$ at $x_k$ is used. This can be phrased as
\begin{equation}
\label{eqn:projected_linesearch}
q_k(P(x_k - \alpha \nabla f(x_k))), \quad \text{for} \ \alpha \in (0,\infty)
\end{equation}
where
\begin{equation}
\label{eqn:projector}
P(x)_i = \begin{cases} l_i, &\text{if} \ x_i < l_i, \\ 
x_i, &\text{if} \ l_i \le x_i \le u_i, \\
u_i, &\text{if} \ u_i < x_i. \end{cases}.
\end{equation}
The latter function is a projector into the box-constraints.

The algorithm just choose the first local minimiser along this line \eqref{eqn:projected_linesearch}, called \textit{\textsc{Cauchy} point}.

\begin{definition}[\textsc{Cauchy} point]
	Let $g: [0,\infty) \rightarrow \R$. Then we call a local minimiser $x \in \R$ of $g$ the \textit{Cauchy point of $g$}, denoted $x_c(g)$, if for every local minimiser $y \in \R \setminus \set{x}$ of $g$ it holds that $x < y$.
\end{definition}

It is then possible, to find the \textsc{Cauchy} point $x_c$ of the projected line search \eqref{eqn:projected_linesearch} analytically. Some of the components of $x_c$ might reach their bounds during this process. We keep these components fixed while doing an unconstrained minimisation of the quadratic model $q_k$ starting at $x_c$ using the remaining components. It is not needed to be exact during this minimisation. Projection as per \eqref{eqn:projector} is then used to fix violated bounds yielding a value $\tilde{x}_{k+1}$.

Eventually, the new descent direction for a line search is given by $p_k:= \tilde{x}_{k+1}- x_k$. For more details see \cite[ch. 18.7]{numericaloptmization}.

%% file: figures/BFGS/scheme.tex
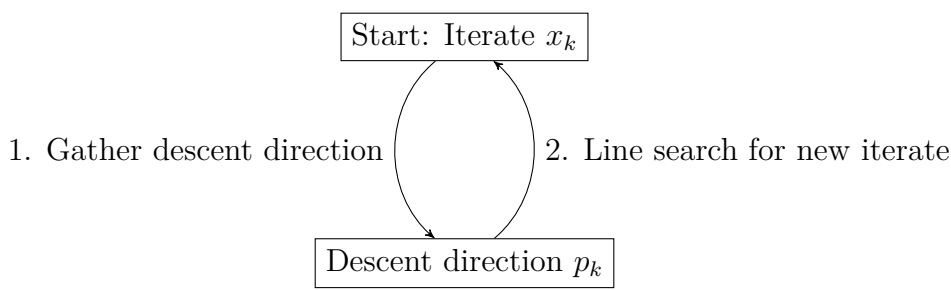
\begin{figure}[b]
	\centering			
	\tikzsetnextfilename{optimisation_scheme}
	\begin{center}
		\begin{tikzpicture}[node distance=3cm, ->, >=stealth']
		\node[draw](source) {Start: Iterate $x_k$};
		\node[draw] (IR) [below of=source] {Descent direction $p_k$};
		
		\path[->]	
		(source) edge[bend right=50] node[left] {1. Gather descent direction} (IR)	
		(IR) edge[bend right=50] node[right] {2. Line search for new iterate} (source);
		\end{tikzpicture}
	\end{center}	
	\caption[General optimistion]{General scheme of minimisation}
	\label{fig:optimisation_scheme}
\end{figure}

%% file: figures/BFGS/contour.tex
\begin{figure}
	\centering			
	\includegraphics{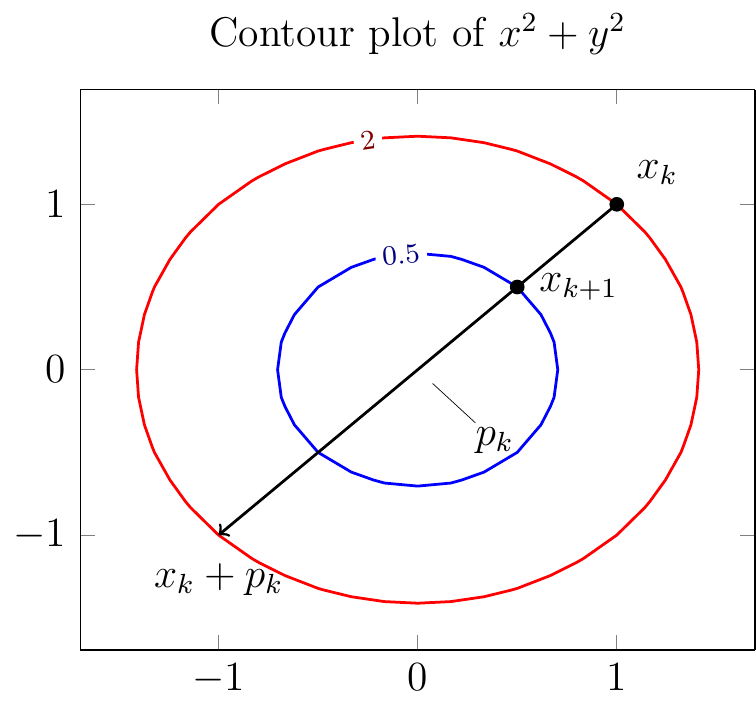}
	\caption[Example of optimisation]{Steps of optimisation in Example~\ref{exmpl:loss}.}
	\label{fig:contour}
\end{figure}

%% file: figures/BFGS/sufficient_decrease.tex
\begin{figure}
	\centering			
	\tikzsetnextfilename{sufficient_decrease}
	\begin{tikzpicture}
	\begin{axis}[
		xlabel={$\alpha$}, 
		ymajorticks=false,
		xmin=0,
		xmax=2.75,
		domain=0:2.75,
		xtick={0,0.5,1,1.5,2,2.5},
		samples=200,
		x axis line style={name path=axis},
		width=1.\textwidth,
		height=0.3\textheight]
		\addplot [color=GAUSSC, thick] {10. - 28.3333*x + 50.75*x^2 - 31.5833*x^3 + 7.*x^4 - 0.333333*x^5} node [pos=0.,pin={-45:\color{black}$f(x_k)$},inner sep=0pt] {};
		\addplot [color=ZOLOC, thick, dashed] {10 - 2.5*x};	
		\addlegendentry{$f(x_k + \alpha p_k) = f(x_{k+1})$}
		\addlegendentry{$c_1 \alpha \nabla f(x_k)^T p_k$}
	\end{axis}
	\end{tikzpicture}	
	\caption[Sufficient decrease]{Illustration of the first \textsc{Wolfe} condition, called \enquote{sufficient decrease}.}
	\label{fig:sufficient_decrease}
\end{figure}
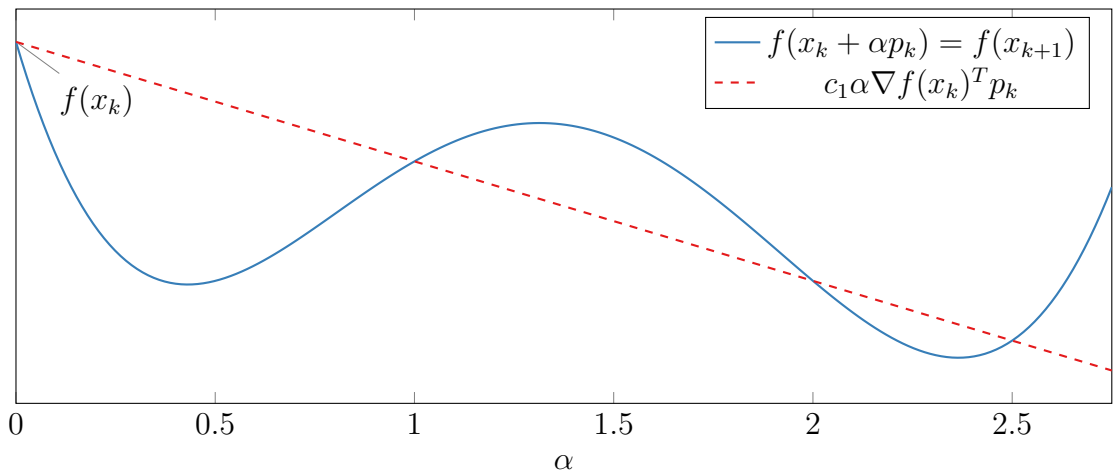

%% file: figures/BFGS/BFGS_scheme.tex
\begin{figure}[b]
	\centering			
	\tikzsetnextfilename{BFGS_scheme}
	\begin{center}
		\begin{tikzpicture}[node distance=1cm, ->, >=stealth']
		\node (input) {Input: $B_k$, current iterate $x_k$};
		\node [draw] (compute) [below of=input] {\BFGS{} model};
		\node (output) [below of=compute] {Output: Descent direction $p_k$};
		
		\path[->]	
		(input) edge (compute)	
		(compute) edge (output);
		\end{tikzpicture}
	\end{center}	
	\caption[Scheme BFGS Algorithm]{Scheme of the \BFGS{} model.}
	\label{fig:BFGS_scheme}
\end{figure}

%% file: figures/BFGS/BFGS_scheme_final.tex
\begin{figure}[b]
	\centering			
	\tikzsetnextfilename{BFGS_scheme_final}
	\begin{center}
		\begin{tikzpicture}[node distance=4cm, ->, >=stealth']
		\node [draw] (compute) {\BFGS{} model};
		\node [draw] (compute2) [below left of=compute] {\BFGS{} Update};
		\node [draw] (linesearch) [below right of=compute] {Line search};
		
		\path[->]	
		(compute) edge node[right] {Descent direction $p_k$} (linesearch)
		(linesearch) edge node[above] {Next iterate $x_{k+1}$} (compute2)
		(compute2) edge node[left] {$B_k$, current iterate $x_k$} (compute);
		\end{tikzpicture}
	\end{center}	
	\caption[BFGS optimisation]{Scheme of \BFGS{} minimisation.}
	\label{fig:BFGS_scheme_final}
\end{figure}

%% file: figures/BFGS/box.tex
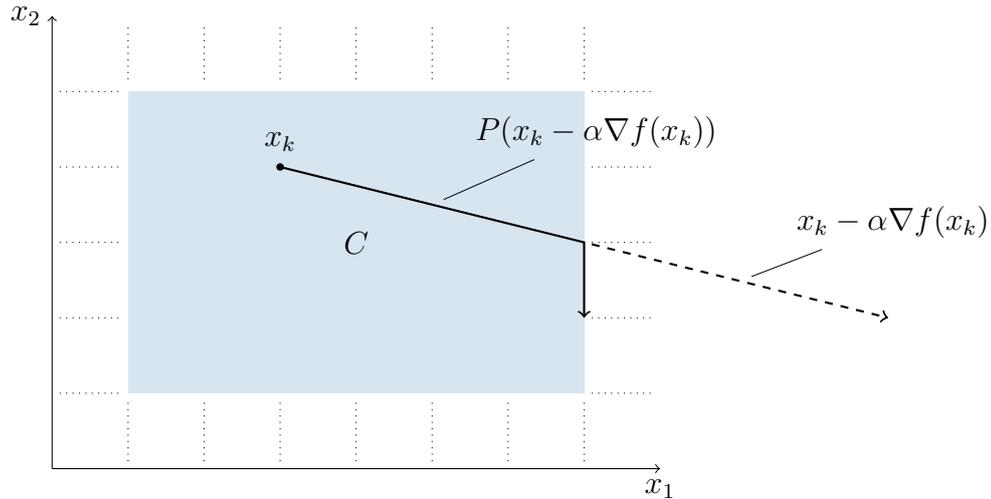
\begin{figure}
	\centering			
	\tikzsetnextfilename{box}
	\begin{tikzpicture}
		[
			grid/.style={black,dotted},
			axis/.style={->,black},
			descent/.style={thick,black,->},
			dot/.style 2 args={circle,inner sep=1pt,fill,label={#2:#1},name=#1}
		]
		\draw[axis] (-4,-3) -- (4,-3) node[anchor=north]{$x_1$};
		\draw[axis] (-4,-3) -- (-4,3) node[anchor=east]{$x_2$};
		\draw[grid] (-3.9,-2.9) grid[xstep=1, ystep=1]  (3.9,2.9);		
	
		\fill[fill=white] (-3.1,-2.1) rectangle (3.1,2.1); 
		\fill[fill=GAUSSC!20!white] (-3,-2) rectangle (3,2) node[pos=.5] {$C$};
		
		\coordinate (X) at (-1, 1); 
		\coordinate (B) at (3, 0); 
		\coordinate (E1) at (3,-1); 
		\coordinate (E2) at (7,-1); 
		\draw[descent] (X) node [dot={$x_k$}{above}] {} -- 
					   node [pin={[pin edge={-,black}]55:$P(x_k- \alpha \nabla f(x_k))$}] {} (B) --
					   (E1);		
		\draw[descent,dashed,thick] ($ (B)!0.025!(E2) $) -- node [pin={[pin edge={-,solid,black}]35:$x_k - \alpha \nabla f(x_k)$}] {} (E2);
	\end{tikzpicture}
	\caption{An example of a box-constraints given by a set $C$ and the projection ${P(x_k - \alpha \nabla f(x_k))}$.}
	\label{fig:box}
\end{figure}

%% file: figures/BFGS/scheme_constraint.tex
\begin{figure}[b]
	\centering			
	\tikzsetnextfilename{BFGS_scheme_contraint}
	\begin{center}
		\begin{tikzpicture}[node distance=4cm, ->, >=stealth']
		\node [draw] (compute) {\textbf{Constrained} \BFGS{} model};
		\node [draw] (compute2) [below left of=compute] {\BFGS{} Update};
		\node [draw] (linesearch) [below right of=compute] {Line search};
		
		\path[->]	
		(compute) edge node[right] {Descent direction $p_k$} (linesearch)
		(linesearch) edge node[above] {Next iterate $x_{k+1}$} (compute2)
		(compute2) edge node[left] {Current iterate $x_k$} (compute);
		\end{tikzpicture}
	\end{center}	
	\caption[Constrained BFGS optimisation]{Scheme of constrained \BFGS{} minimisation.}
	\label{fig:BFGS_scheme_contraint}
\end{figure}
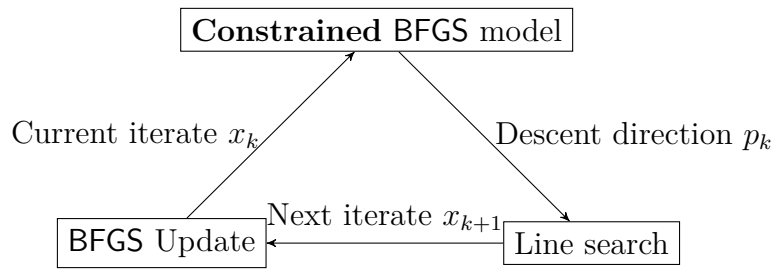